\documentclass[11pt]{article}

\usepackage{amsmath, amsfonts, amssymb,latexsym,wasysym,graphicx} \input epsf.sty

\newtheorem{Theorem}{Theorem}[section] \newtheorem{Lemma}{Lemma}[section]
\newtheorem{Proposition}[Lemma]{Proposition} \newtheorem{Corollary}[Lemma]{Corollary}
 \newtheorem{Definition}[Lemma]{Definition}
\newtheorem{Remark}[Lemma]{Remark} 

\newcommand{\BEQ}{\begin{equation}}     
\newcommand{\BEA}{\begin{eqnarray}}
\newcommand{\BD}{\begin{displaymath}} \newcommand{\EEQ}{\end{equation}}       
\newcommand{\EEA}{\end{eqnarray}} \newcommand{\ED}{\end{displaymath}} 
\newcommand{\del}{\delta} \newcommand{\Del}{\Delta}
\newcommand{\eps}{\varepsilon}          

\newcommand{\R}{\mathbb{R}}  \newcommand{\Z}{\mathbb{Z}}
\newcommand{\N}{\mathbb{N}}

\newcommand{\SkI}{{\mathrm{SkI}}} \newcommand{\supp}{{\mathrm{supp}}} \renewcommand{\H}{{\bf H}}
\newcommand{\Roo}{{\mathrm{Roo}}} \newcommand{\Lea}{{\mathrm{Lea}}}

 \def\esper{{\mathbb{E}}} \def\T{{\mathbb{T}}} \def\F{{\mathbb{F}}}
\def\Var{{\mathrm{Var}}}

\def\shuffle{\pitchfork}

\newcommand{\tdun}[1]{\begin{picture}(10,5)(-2,-1) \put(0,0){\circle*{2}} \put(3,-2){\tiny #1}
\end{picture}}

\newcommand{\tddeux}[2]{\begin{picture}(12,5)(0,-1) \put(3,0){\circle*{2}} \put(3,0){\line(0,1){5}}
\put(3,5){\circle*{2}} \put(6,-2){\tiny #1} \put(6,3){\tiny #2} \end{picture}}

\newcommand{\tdtroisun}[3]{\begin{picture}(20,12)(-5,-1) \put(3,0){\circle*{2}} \put(-0.65,0){$\vee$}
\put(6,7){\circle*{2}} \put(0,7){\circle*{2}} \put(5,-2){\tiny #1} \put(9,5){\tiny #2} \put(-5,5){\tiny
#3} \end{picture}} \newcommand{\tdtroisdeux}[3]{\begin{picture}(12,12)(-2,-1) \put(0,0){\circle*{2}}
\put(0,0){\line(0,1){5}} \put(0,5){\circle*{2}} \put(0,5){\line(0,1){5}} \put(0,10){\circle*{2}}
\put(3,-2){\tiny #1} \put(3,3){\tiny #2} \put(3,9){\tiny #3} \end{picture}}

\newcommand{\tdquatredeux}[4]{\begin{picture}(20,20)(-5,-1) \put(3,0){\circle*{2}} \put(-.65,0){$\vee$}
\put(6,7){\circle*{2}} \put(0,7){\circle*{2}} \put(0,14){\circle*{2}} \put(0,7){\line(0,1){7}}
\put(5,-2){\tiny #1} \put(9,5){\tiny #2} \put(-5,5){\tiny #3} \put(-5,12){\tiny #4} \end{picture}}

\newcommand{\eop}{\hfill $\Box$}        

\newcommand{\II}{{\rm i}}               
\newcommand{\half}{{1\over 2}}          

\catcode`\@=11 \def\numberbysection{\@addtoreset{equation}{section}
        \def\theequation{\thesection.\arabic{equation}}}
\numberbysection

\begin{document}

\vspace*{1.5cm} \begin{center} {\Large \bf A renormalized  rough path over fractional Brownian motion}
\end{center}

\vspace{2mm} \begin{center} {\bf  J\'er\'emie Unterberger} \end{center}

\vspace{2mm} \begin{quote}

\renewcommand{\baselinestretch}{1.0} \footnotesize {We construct in this article a  rough path over
fractional Brownian motion with arbitrary Hurst index  by (i)
 using the
Fourier normal ordering algorithm introduced in \cite{Unt-Holder} to reduce the problem to that of
regularizing tree iterated integrals and (ii) applying  the Bogolioubov-Parasiuk-Hepp-Zimmermann (BPHZ)
renormalization algorithm to Feynman diagrams representing tree iterated integrals.
 }
\end{quote}

\vspace{4mm} \noindent {\bf Keywords:} fractional Brownian motion, rough paths, H\"older continuity,
renormalization,  Hopf algebra of decorated rooted trees, shuffle algebra


\tableofcontents

\newpage


\section{Introduction}


Consider a $d$-dimensional continuous path $t\mapsto \Gamma_t=(\Gamma_t(1),\ldots,\Gamma_t(d))$,
$t\in\R$. Assume $\Gamma$ is {\em not} differentiable, but only $\alpha$-H\"older for some
$\alpha\in(0,1)$.  Rough path theory answers positively the following related two questions, in
particular: (i) can one integrate a (sufficiently regular) one-form along $\Gamma$ ? (ii) can one solve
differential equations driven by $\Gamma$ ? The solution of these relies on the definition of a
so-called {\em rough path over $\Gamma$}, denoted by ${\bf
\Gamma}=({\bf\Gamma}^{ts}(i_1,\ldots,i_n))_{1\le n\le \lfloor 1/\alpha\rfloor,1\le i_1,\ldots,i_n\le
d}$, which is a {\em substitute for iterated integrals} $\int_s^t
d\Gamma_{t_1}(i_1)\ldots\int_s^{t_{n-1}} d\Gamma_{t_n}(i_n)$ for $n=1,\ldots,\lfloor 1/\alpha\rfloor$,
defined both by {\em algebraic} and {\em regularity} properties (see Definition \ref{def:rough-path}).
Rough path solutions are remarkably well-behaved with respect to controllability and numerical schemes,
and the construction is robust enough to extend to a variety of settings.

Given this, it is important to know how to construct a rough path. A first  answer to this problem has
been given by T. Lyons and N. Victoir \cite{LyoVic07}. However, their construction is non-canonical
(actually, it uses the axiom of choice) and does not provide a closed formula, which bars the way to
applications to non-pathwise results for stochastic processes for instance. Among these, fractional
Brownian motion (fBm for short) is probably the one which has drawn most attention, probably because it
is the simplest non-trivial example. This Gaussian, self-similar processes, depending on a regularity
index $\alpha\in (0,1)$ called Hurst index, has $\alpha^-$-H\"older (i.e. $(\alpha-\eps)$-H\"older for
every $\eps>0$) paths. Consider a $d$-dimensional fBm, $B_t=(B_t(1),\ldots,B_t(d))$, with $d\ge 2$ (the
one-dimensional case is much simpler and has been solved earlier \cite{GNRV}). Classical results  imply
that the natural iterated integrals of the piecewise linear \cite{CQ02} or analytic \cite{Unt-1}
approximation of fBm converge to a rough path over $B$ if and only if $\alpha>1/4$. The search for other
Gaussian approximations with converging iterated integrals has failed up to now, and recent
investigations have turned (i) either to {\em non-Gaussian approximations}, using the tools of {\em
constructive quantum field theory} \cite{MagUnt} (including renormalization); or (ii) to {\em
"algebraic" rough paths}, i.e. substitute for iterated integrals in the above sense, satisfying the
required algebraic and regularity properties, but not given by any explicit approximation \footnote{Such
approximations -- using pieces of sub-riemannian geodesics --  have been shown to exist in general, but
are not very explicit \cite{FV}}. It is the second approach that we pursue in this article, but always
keeping an eye on the first one, as we shall see.

This approach relies on a combinatorial algorithm called {\em Fourier normal ordering}. Initially, it
was conceived as a splitting into sectors of the domain of integration in Fourier coordinates which
produces naturally H\"older bounds \cite{Unt09ter}.
 For iterated integrals of
lowest orders at least, it appeared clearly that recombining regularized iterated integrals defined
within each sector gave a quantity satisfying the algebraic properties required for a rough path. With
the time, it became clear that Fourier normal ordering made it possible to separate the rough path
construction problem into two questions of a totally different nature:

-- the first one consists in {\em regularizing tree iterated integrals} or more precisely {\em tree
skeleton integrals} -- restricted to the above Fourier sectors --, which are natural combinatorial
extensions of iterated integrals indexed by decorated trees;

-- the second one consists in showing that one may reconstruct in a canonical way a rough path out of
these data.

It turns out that rough path construction is a {\em very undetermined problem}, since in some sense {\em
any} regularization scheme (including the brutal-force regularization by zero, except for first-order
integrals) gives in the end a {\em formal rough path}, i.e. a set of  quantities satisfying the algebraic requirements. It seems also rather clear
-- without pretending to make this a formal statement -- that regularized tree skeleton integrals with
the correct H\"older regularity should yield by recombination a rough path with the correct H\"older
regularity.

\bigskip

Taking for granted the combinatorial part of  Fourier normal ordering -- which we briefly recall in
section 1 for completeness -- one is naturally led to decide {\em which} regularization scheme is most
natural. We belive that the only possible answer to
 this question is to provide a natural approximation scheme leading to the corresponding
rough path, which leads us back to the first approach -- still under way -- using quantum field theory
methods \cite{MagUnt}. Its perturbative formulation is based on the Bogolioubov-Parasiuk-Hepp-Zimmermann
(BPHZ for short) renormalization scheme for Feynman diagrams \cite{Hepp}. To say things shortly, this is a recursive
method to discard nested divergences, depending on the choice of a regularization scheme for diagrams
without sub-divergences. Usually, the renormalization is implemented by a change of the parameters of
the measure. Here, however, the theory is a priori free, i.e. Gaussian, and such an implementation is
impossible without changing the definition of the underlying process, see again \cite{MagUnt} for a way
out of this. Hence any Gaussian renormalization is in some sense arbitrary. Nevertheless it seems
natural to mimic the renormalization schemes of quantum field theory in the following way. The {\em
variance of iterated integrals} may be represented as {\em Feynman diagrams}; iterated integrals
themselves are represented by {\em Feynman "half-diagrams"} and evaluated by integrating some {\em
deterministic  kernel} against a multi-dimensional Brownian motion. Renormalizing directly Feynman
diagrams, as mentioned above, leads us  to the non-Gaussian constructive field theory approach. Instead,
we choose here to renormalize the {\em kernel}, still by the same BPHZ algorithm, which is a
non-conventional approach. This yields directly a renormalized {\em random variable} in the same chaos
as the original, unrenormalized quantity, which is proved to enjoy the required H\"older regularity.

\bigskip

Our main result may be stated as follows.

\begin{Theorem} \label{th:0.1} Let $\alpha\in(0,1)$ such that $1/\alpha\not\in\N$. Let ${\bf
B}^{ts}(i_1,\ldots,i_n):=J^{ts}_B(i_1,\ldots,i_n)$, $n=1,\ldots,\lfloor 1/\alpha \rfloor$ be the random
variable in the $n$-th chaos of fBm, defined in Proposition \ref{lem:barchi-char} and Definition
\ref{def:3.2}. Then:

\begin{enumerate} \item $||{\bf B}^{ts}(i_1,\ldots,i_n)||_{2,n\alpha}:=\sup_{s,t\in [0,T]} \frac{ |{\bf
B}^{ts}(i_1,\ldots,i_n)| }{|t-s|^{n\alpha}}$ is an $L^2$ random variable. \item ${\bf B}:=({\bf
B}^{ts}(i_1,\ldots,i_n))_{1\le n\le \lfloor 1/\alpha\rfloor, 1\le i_1,\ldots,i_n\le d}$ satisfies the
Chen and shuffle properties
   (see Definition \ref{def:rough-path}).
\end{enumerate}

Hence ${\bf B}$ is an $\alpha^-$-H\"older rough path over $B$.

\end{Theorem}

{\bf Remarks.}

\begin{enumerate} \item If $1/\alpha\in\N$, and $\kappa<\alpha$ is chosen as close to $\alpha$ as
desired, then Theorem \ref{th:0.1} applies to $B$ seen as a $\kappa$-H\"older path, and yields a
$\kappa^-$-H\"older rough path over $B$.

\item Property (1) in Theorem \ref{th:0.1} is a consequence of the estimates \BEQ \esper |{\bf
    B}^{ts}(i_1,\ldots,i_n)|^2\le C|t-s|^{2n\alpha},\quad n\le \lfloor 1/\alpha\rfloor \EEQ proved
    in section 5, as follows from the Garsia-Rodemich-Rumsey lemma \cite{Gar} and from the
    equivalence of $L^p$-norms for variables in a fixed Gaussian chaos, see \cite{Unt-fBm}, section
    1, for details.

\end{enumerate}

Here is a plan of the article. We start in section 1 by recalling the fundamentals of the Fourier normal
ordering algorithm, refering to \cite{Unt-Holder,FoiUnt}  for a complete treatment. The correspondence
with Feynman diagrams and half-diagrams is explained in Section 2. The systematics of renormalization,
including its multi-scale version which has been acknowledged as the quickest way to get estimates, is
recalled in Section 3. We use a classical multi-scale expansion to derive a general bound for Feynman
diagrams in  Section 4. We conclude in Section 5 by proving the H\"older estimates for the rough path
and adding some remarks on related previous attempts and on possible extensions to  general H\"older
paths.


\section{The  Fourier normal ordering algorithm}


Let $\Gamma=(\Gamma_t(1),\ldots,\Gamma_t(d)):\R\to\R^d$ be some continuous path, compactly supported in
$[0,T]$. Assume that $\Gamma$ is not differentiable, but only $\alpha$-H\"older for some $0<\alpha<1$,
i.e. bounded in the ${\mathcal C}^{\alpha}$-norm, \BEQ ||\gamma||_{{\mathcal C}^{\alpha}}:=\sup_{t\in
[0,T]} ||\Gamma_t||+\sup_{s,t\in[0,T]} \frac{||\Gamma_t -\Gamma_s||}{|t-s|^{\alpha}}.\EEQ Then iterated
integrals of $\Gamma$  are not canonically defined. As explained in the Introduction, rough path theory
may be seen as a black box taking as input some lift of $\Gamma$ called {\em rough path over} $\Gamma$,
 producing e.g. solutions of differential equations driven by $\Gamma$.


\subsection{Rough paths and iterated integrals}


 The usual
definition of a rough path is the following. We let in the sequel
 $\lfloor 1/\alpha\rfloor$ be the entire part of $1/\alpha$.

\begin{Definition} \label{def:rough-path}

A rough path over $\Gamma$ is a functional $J_{\Gamma}^{ts}(i_1,\ldots,i_n)$, $n\le \lfloor
1/\alpha\rfloor$, $i_1,\ldots,i_n\in\{1,\ldots,d\}$, such that $J_{\Gamma}^{ts}(i)=
\Gamma_t(i)-\Gamma_s(i)$ are the increments of $\Gamma$, and the following 3 properties are satisfied:

\begin{itemize} \item[(i)] (H\"older continuity) $J_{\Gamma}^{ts}(i_1,\ldots,i_n)$ is $n\alpha$-H\"older
continuous as a function of two variables, namely, $\sup_{s,t\in\R}
\frac{|J_{\Gamma}^{ts}(i_1,\ldots,i_n)|}{|t-s|^{\alpha}}<\infty.$ \item[(ii)] (Chen property) \BEQ
J_{\Gamma}^{ts}(i_1,\ldots,i_n)=J_{\Gamma}^{tu}(i_1,\ldots,i_n)+ J_{\Gamma}^{us}(i_1,\ldots,i_n)+
\sum_{n_1+n_2=n} J_{\Gamma}^{tu}(i_1,\ldots,i_{n_1})J_{\Gamma}^{us}(i_{n_1+1},\ldots,i_{n});\EEQ

\item[(iii)] (shuffle property) \BEQ
    J_{\Gamma}^{ts}(i_1,\ldots,i_{n_1})J_{\Gamma}^{ts}(j_1,\ldots,j_{n_2})=\sum_{\vec{k}\in
    Sh(\vec{i},\vec{j})} J_{\Gamma}^{ts}(k_1,\ldots,k_{n_1+n_2}),\EEQ where $Sh(\vec{i},\vec{j})$ --
    the set of shuffles of the words $\vec{i}$ and $\vec{j}$ -- is the subset of permutations of the
    union of the lists $\vec{i},\vec{j}$ leaving unchanged the order of the sublists $\vec{i}$ and
    $\vec{j}$. For instance, $J^{ts}_{\Gamma}(i_1,i_2)J^{ts}_{\Gamma}(j_1)=
    J^{ts}_{\Gamma}(i_1,i_2,j_1)+J^{ts}_{\Gamma}(i_1,j_1,i_2)+J^{ts}_{\Gamma}(j_1,i_1,i_2).$
    \end{itemize}
    
    A {\em formal} rough path over $\Gamma$ is a functional satisfying all the above properties
    {\em except} H\"older continuity (i).

\end{Definition}

In particular, if $\Gamma$ is smooth, then  its natural iterated integrals \BEQ
I^{ts}_{\Gamma}(i_1,\ldots,i_n):=\int_s^t d\Gamma_{t_1}(i_1)\ldots\int_s^{t_{n-1}} d\Gamma_{t_n}(i_n)
\EEQ satisfy properties (ii) and (iii).

\bigskip

These two algebraic axioms  may be rewritten in a Hopf algebraic language. Let us say a few words about
it. The reader who is allergic to algebra may just read Definition  \ref{def:H} and Proposition
\ref{prop:permutation-graph}, skip the rest of the section and jump to the end of subsection 1.4.
 However,  this language has proved to be very useful both from a theoretic
and a practical point of view \cite{Unt-Holder,FoiUnt}.

\begin{Definition}[Hopf algebra of decorated rooted trees] \label{def:H}

\begin{itemize} \item[(i)] A decorated rooted tree is a tree with a distinguished vertex called root
(drawn growing up from the root to the top), provided with a decoration for each vertex. In this
article, decorations are always assumed to range in the set $\{1,\ldots,d\}$. The set of trees is
denoted by $\cal T$. The commutative product $\T_1.\T_2$ of two trees yields the  forest with the two
connected components $\T_1$ and $\T_2$. The algebra over $\R$ generated by trees is denoted by ${\bf
H}$, and the linear subspace of forests with $n$ vertices by ${\bf H}(n)$.

\item[(ii)] If $w$ is a descendant of $v$ (i.e. $w$ is above $v$) then one writes
    $w\twoheadrightarrow v$. One says that $v$ is connected to $w$ (a symmetric relation) if either
    $w=v$, $w\twoheadrightarrow v$ or $v\twoheadrightarrow w$. A subset of vertices $\vec{v}\subset
    V(\T)$ is an {\em admissible cut} if $(v,w\in \vec{v}, v\not=w)\Rightarrow (v$ is not connected
    to $w)$. If $\vec{v}$ is admissible, which we write $\vec{v}\models V(\T)$, then
    $Roo_{\vec{v}}\T$ is the subforest with vertices $\{w\in V(\T); \exists v\in \vec{v},
    v\twoheadrightarrow w\}$, while $Lea_{\vec{v}}\T$ is the subforest with the complementary set of
    vertices. Note that $Roo_{\vec{v}}\T$ is a tree if $\T$ is a tree.

\item[(iii)] Define \BEQ \Del(\T)=\sum_{v\models V(\T)} Roo_{\vec{v}}\T \otimes Lea_{\vec{v}}\T.\EEQ
    Then $\bf H$ equipped with $\Del:{\bf H}\to {\bf H}\otimes {\bf H}$ is a coproduct. For
    instance, \BEQ \Del(\tdtroisun{a}{c}{b})=\tdtroisun{a}{c}{b}\otimes 1+1\otimes
    \tdtroisun{a}{c}{b}+ \tddeux{a}{b}\otimes \tdun{c} +\tddeux{a}{c}\otimes
    \tdun{b}+\tdun{a}\otimes \tdun{b}\tdun{c}\EEQ

\item[(iv)] $\H$ has an antipode $\bar{S}$, defined inductively by \BEQ \bar{S}(1)=1,\quad
    \bar{S}(\T)=-\T-\sum_{\vec{v}\models V(\T),\vec{v}\not=\emptyset} Roo_{\vec{v}}\T \cdot
    \bar{S}(Lea_{\vec{v}}\T).\EEQ

\end{itemize} \end{Definition}

We shall also need the following Hopf algebra in order to encode the shuffle property.

\begin{Definition}[shuffle algebra]

\begin{itemize} \item[(i)] Let ${\bf Sh}$ be the shuffle algebra with decorations in $\{1,\ldots,d\}$,
i.e. the set of words $(i_1 \ldots i_n)$, $i_1,\ldots,i_n\in\{1,\ldots,d\}$, with product \BEQ
(i_1\ldots i_{n_1})\shuffle (j_1\ldots j_{n_2})=\sum_{\vec{k}\in Sh(\vec{i},\vec{j})} (k_1\ldots
k_{n_1+n_2}).\EEQ An element of $\bf Sh$ is naturally represented as a trunk tree decorated by
$\ell=(\ell(1),\ldots,\ell(n))$ from the root to the top. For instance $(i_1 i_2 i_3)=
\tdtroisdeux{$i_1$}{$i_2$}{$i_3$}$ is decorated by $\ell(j)=i_j$, $j=1,2,3$. \item[(ii)] $\bf Sh$
equipped with the restriction of the coproduct $\Del$ of $\bf H$ to trunk trees, and with the antipode
$S((i_1\ldots i_n))=-(i_n\ldots i_1)$, is a Hopf algebra. It holds: $\Del((i_1\ldots i_n))=\sum_{k=0}^n
(i_1\ldots i_k)\otimes (i_{k+1}\ldots i_n).$

\end{itemize} \end{Definition}

One of the links between these two algebras is given by the following Proposition.

\begin{Proposition}[projection morphism]

Let $\theta:{\bf H}\to {\bf Sh}$ be the projection Hopf morphism given by associating to a tree $\T$ the
sum of the trunk trees  $\mathfrak{t}$ with same decorations such that \BEQ (v\twoheadrightarrow w\
{\mathrm{in}} \ \T)\Rightarrow (v\twoheadrightarrow w\ {\mathrm{in}}\ \mathfrak{t}).\EEQ

For instance
$\theta(\tdtroisun{$a$}{$c$}{$b$})=\tdtroisdeux{$a$}{$b$}{$c$}+\tdtroisdeux{$a$}{$c$}{$b$}.$

\end{Proposition}

 Indexing
the $J_{\Gamma}^{ts}(i_1,\ldots,i_n)$ by trunk trees $\T\in {\bf Sh}$ with decoration $\ell(j)=i_j$,
$j=1,\ldots,n$, properties (ii) and (iii) in Definition \ref{def:rough-path} are equivalent to
\begin{itemize} \item[{\em (ii)bis}] \BEQ J^{ts}_{\Gamma}(\T)= \sum_{\vec{v}\models V(\T)}
J^{tu}_{\Gamma}(\Roo_{\vec{v}}(\T)) J^{us}_{\Gamma}(\Lea_{\vec{v}}(\T)),\quad \T\in{\bf Sh};\EEQ in
other words, $J^{ts}_{\Gamma}=J^{tu}_{\Gamma} \ast J^{us}_{\Gamma}$ for the shuffle convolution defined
in subsection 5.2; \item[{\em (iii)bis}] \BEQ J^{ts}_{\Gamma}(\T)J^{ts}_{\Gamma}(\T')=J^{ts}_{\Gamma}(\T
\shuffle \T'), \quad \T,\T'\in {\bf Sh}.\EEQ In other words, $J^{ts}_{\Gamma}$ is a character of ${\bf
Sh}$. \end{itemize}

Such a functional indexed by {\em trunk trees} extends easily to a general {\em tree-indexed functional}
or {\em tree-indexed rough path} by setting $\bar{J}^{ts}_{\Gamma}(\mathbb{T}):=
J^{ts}_{\Gamma}\circ\theta(\mathbb{T})$.  Since $\theta$ is a Hopf algebra morphism, one gets
immediately the generalized properties \begin{itemize} \item[{\em (ii)ter}]
$\bar{J}^{ts}_{\Gamma}=\bar{J}^{tu}_{\Gamma} \ast \bar{J}^{us}_{\Gamma}$ for the
 convolution of $\H$, i.e.
\BEQ \bar{J}^{ts}_{\Gamma}(\T)=\sum_{\vec{v}\models V(\T)} \bar{J}^{tu}_{\Gamma}(\Roo_{\vec{v}}\T)
\bar{J}^{us}_{\Gamma}(\Lea_{\vec{v}}\T) ;\EEQ \item[{\em (iii)ter}]
$\bar{J}^{ts}_{\Gamma}(\T)\bar{J}^{ts}_{\Gamma}(\T')= \bar{J}^{ts}_{\Gamma}(\T.\T')$, in other words,
$\bar{J}^{ts}_{\Gamma}$ is a character of $\H$. \end{itemize}

Properties (ii), (iii) and their generalizations are satisfied for the usual integration operators
$I^{ts}_{\Gamma}$ and their tree extension $\bar{I}^{ts}_{\Gamma}$, provided $\Gamma$ is a smooth path
so that iterated integrals make sense \cite{Gu}.

Let us give an explicit formula for tree iterated integrals. Let $\T$ be e.g. a tree, and index its
vertices as $1,\ldots,n$, so that $(i\twoheadrightarrow j)\Rightarrow (i>j)$.
 Denoting by $i^-$ the ancestor of the vertex $i$ in $\T$, one has
\BEQ \bar{I}^{ts}_\Gamma(\T)=\int_s^t d\Gamma_{x_1}(\ell(1))\int_s^{x_{2^-}} d\Gamma_{x_2}(\ell(2))
\ldots \int_s^{x_{n^-}} d\Gamma_{x_n}(\ell(n)). \label{eq:barIts} \EEQ

\begin{Remark} \label{rem:1.5} Note that (\ref{eq:barIts}) obviously does not depend on the choice of
the vertex indexation. We call this {\em invariance under indexation of the vertices}, or {\em
naturality property}. We may rephrase it saying that $\bar{I}^{ts}_{\Gamma}(\T)$ depends only on the
{\em topology} of $\T$. The same property applies to every natural construction and is required in
Definition \ref{def:barchi-char}. \end{Remark}

\bigskip

Suppose now one wishes to construct a rough path over $\Gamma$, and concentrate on the algebraic
properties (ii), (iii) of Definition \ref{def:rough-path}.
 Assume one has constructed characters of ${\bf Sh}$,
$J^{ts_0}_{\Gamma}$, $t\in [0,T]$ with $s_0$ fixed -- in other words, a one-time functional satisfying
the usual shuffle property (iii) --,
 such that $J_{\Gamma}^{ts_0}(i)=\Gamma_t(i)-\Gamma_{s_0}(i)$,
 then one immediately checks that
$J^{ts}_{\Gamma}:=J^{ts_0}_{\Gamma}\ast (J^{ss_0}_{\Gamma}\circ \bar{S})$ satisfies properties (ii)bis
and (iii)bis. Namely (by the definition of the antipode) $J^{ts}_{\Gamma}:=J^{ts_0}_{\Gamma}\ast
(J^{ss_0}_{\Gamma}\circ \bar{S})$ is equivalent to the Chen property $J^{ts}_{\Gamma}\ast
J^{ss_0}_{\Gamma} =J^{ts_0}_{\Gamma}$.  So the only difficult part consists in defining some regularized
character of ${\bf Sh}$ satisfying the regularity properties (i).


\subsection{Fourier transform and skeleton integrals}


Instead of regularizing iterated integrals, $I_{\Gamma}^{ts_0}\rightsquigarrow J_{\Gamma}^{ts_0}$ with
$s_0$ fixed, we choose to regularize {\em skeleton integrals}, $\SkI_{\Gamma}^t$, which are analogues of
iterated integrals but depending naturally on a single argument, defined by using  Fourier transform.

\begin{Definition}[skeleton integral] Let \BEA &&  \SkI^t_{\Gamma}(a_1\ldots a_n):= \nonumber\\ &&
\qquad (2\pi)^{-n/2} \int_{\R^n} \prod_{j=1}^n {\mathcal F}\Gamma'_{\xi_j}(a_j)d\xi_j\ \cdot\ \int^t
dx_1\int^{x_1}dx_2\ldots \int^{x_{n-1}} dx_n e^{\II(x_1\xi_1+\ldots+x_n\xi_n)}, \nonumber\\ \EEA where,
by definition, $\int^x e^{\II y\xi} dy=\frac{e^{\II x\xi}}{\II\xi}.$ It may be checked that
$\SkI^t_{\Gamma}$ is a character of ${\bf Sh}$ -- or, in other words, satisfies the shuffle property --,
just as for usual iterated integrals. \end{Definition}

The projection $\theta$ yields immediately a generalization of this notion to tree skeleton integrals,
compare with eq. (\ref{eq:barIts}), \BEA
\overline{\SkI}^t_{\Gamma}(\mathbb{T})=\SkI^t_{\Gamma}\circ\theta(\mathbb{T})=
 \int^t d\Gamma_{x_1}(\ell(1))\int_s^{x_{2^-}} d\Gamma_{x_2}(\ell(2))
\ldots \int^{x_{n^-}} d\Gamma_{x_n}(\ell(n)). \nonumber\\ \EEA

An explicit computation yields (\cite{Unt-Holder}, Lemma 4.5): \BEQ
\overline{\SkI}^t_{\Gamma}(\mathbb{T})=(2\pi)^{-n/2} \II^{-n} \int_{\R^n} \prod_{j=1}^n {\mathcal
F}\Gamma'_{\xi_j}(\ell(j))d\xi_j\ \cdot\ \frac{e^{\II t(\xi_1+\ldots+\xi_n)}}{\prod_{i=1}^n
[\xi_i+\sum_{j\twoheadrightarrow i} \xi_j]}. \label{eq:1.16} \EEQ


\subsection{Fourier normal ordering for smooth paths}


We begin by the following

\begin{Definition}[Fourier projections and measure-splitting] \label{def:splitting}

\begin{itemize}

\item[(i)]

 Let $\mu$ be some signed measure with compact support,
typically, $\mu=\mu_{(\Gamma,\ell)}(dx_1,\ldots,dx_n)=\otimes_{j=1}^n d\Gamma_{x_j}(\ell(j)).$
  Then \BEQ \mu=\sum_{\sigma\in\Sigma_n} \mu^{\sigma}\circ\sigma^{-1},\EEQ
where \BEQ {\cal P}^{\sigma}:\mu\mapsto {\cal F}^{-1}\left( {\bf 1}_{|\xi_{\sigma(1)}|\le \ldots\le
|\xi_{\sigma(n)}|} {\cal F}\mu(\xi_1,\ldots,\xi_n) \right) \EEQ is a Fourier projection, and
$\mu^{\sigma}$ is defined by \BEQ \mu^{\sigma}:={\cal P}^{{\mathrm{Id}}}(\mu\circ\sigma)=({\cal
P}^{\sigma}\mu)\circ \sigma.\EEQ

The set of all measures whose Fourier transform is supported in $\{(\xi_1,\ldots,\xi_n);
|\xi_{1}|\le \ldots\le |\xi_{n}| \}$ will be  denoted by ${\cal P}^+ Meas(\R^n)$. Thus $\mu^{\sigma}
\in  {\cal P}^+ Meas(\R^n)$.

\item[(ii)]
 More generally,
if $T$ is a tree, \BEQ {\cal P}^{\T} Meas(\R^n)=\left\{\mu; \vec{\xi}\in {\mathrm{supp}}({\cal
F}\mu)\Rightarrow
 \left( (i\twoheadrightarrow j)\Rightarrow (|\xi_i|>|\xi_j|)  \right) \right\}.\EEQ
\end{itemize} \end{Definition}

This definition applies in particular to the tensor measures $\mu=\mu_{(\Gamma,\ell)}=\otimes_{i=1,
\ldots,n} d\Gamma_{x_i}(\ell(i))=\otimes_{i=1,\ldots,n} \Gamma'_{x_i}(\ell(i))dx_i$ if
$\ell=(\ell(1),\ldots,\ell(n))$ is the decoration of a trunk tree. Note that even though $\mu$ is a {\em
tensor measure}  in this case, the projected measures $\mu^{\sigma}$ are {\em not}. This forces us to
extend the previous definitions of $I^{ts}_{\Gamma}$, $J^{ts}_{\Gamma}$, $\bar{I}^{ts}_{\Gamma}$,
$\bar{J}^{ts}_{\Gamma}$, $\SkI^{ts}_{\Gamma}$, $\overline{\SkI}^{ts}_{\Gamma}$
 to measure-indexed characters. This
is  straightforward. However, one must then trade decorated trees (or forests) for so-called {\em
heap-ordered trees} (or {\em forests}), i.e. trees without decoration but with indexed vertices
$1,\ldots,n$ such that \BEQ (i\twoheadrightarrow j)\Rightarrow (i>j). \label{eq:ho} \EEQ For instance,
\BEQ \bar{I}^{ts}_\mu(\T)=\int_s^t \int_s^{x_{2^-}} \ldots \int_s^{x_{n^-}} d\mu(x_1,\ldots,x_n).
 \EEQ

{\bf Remark.}  Recall from Remark \ref{rem:1.5} that iterated integrals depend only on the topology of
the tree, which means that \BEQ \bar{I}^{ts}_{\mu}(\T)=\bar{I}^{ts}_{\mu\circ\sigma} (\sigma^{-1}.\T)
\EEQ if $\sigma\in\Sigma_n$ is a reindexation of the vertices preserving the topology of $\T$, i.e. such
that \BEQ (i\twoheadrightarrow j\ {\mathrm{in}}\ \T)\Rightarrow  (i\twoheadrightarrow j\ {\mathrm{in}}\
\sigma^{-1}(\T)).\EEQ

\bigskip

To say things shortly, {\em skeleton integrals} are convenient when using Fourier coordinates, since
they avoid awkward boundary terms such as those generated by  usual integrals, $\int_0^x e^{\II
y\xi}dy=\frac{e^{\II x\xi}}{\II \xi}-\frac{1}{\II\xi}$, which create terms with different homogeneity
degree in $\xi$ by iterated integrations. {\em Measure splitting} gives the {\em relative scales} of the
Fourier coordinates; orders of magnitude of the corresponding integrals may be obtained separately in
each sector $|\xi_{\sigma(1)}|\le \ldots\le |\xi_{\sigma(n)}|$. It turns out that these are easiest to
get after a permutation of the integrations (applying Fubini's theorem) such that {\em innermost (or
rightmost)integrals bear highest Fourier frequencies}. This is the essence of {\em Fourier normal
ordering}.

\begin{Proposition}[permutation graph] \label{prop:permutation-graph}

Let $\mathfrak{T}_n\in {\bf Sh}$ be a trunk tree with $n$ vertices, and  $\sigma\in\Sigma_n$ a
permutation of $\{1,\ldots,n\}$.
  Then there exists a unique element $\T^{\sigma}\in {\bf H}$ called
{\em permutation graph} such that \BEQ
I^{ts}_{\Gamma}(\mathfrak{T}_n)=I^{ts}_{\Gamma}(\T^{\sigma}).\EEQ

\end{Proposition}

Let us give an example. Let $\mathfrak{T}_n=\tdtroisdeux{$a_1$}{$a_2$}{$a_3$}$ and
$\sigma:(1,2,3)\to(2,3,1).$ Then

\begin{eqnarray*} I^{ts}_\Gamma(\mathfrak{T}_n) &=&\int_s^t d\Gamma_{a_1}(x_3) \int_s^{x_3}
d\Gamma_{a_2}(x_1) \int_s^{x_1} d\Gamma_{a_3}(x_2)\\ &=&\int_s^t d\Gamma_{a_2}(x_1) \int_s^{x_1}
d\Gamma_{a_3}(x_2) \int_{x_1}^t d\Gamma_{a_1}(x_3)\\ &=&\int_s^t d\Gamma_{a_2}(x_1) \int_s^{x_1}
d\Gamma_{a_3}(x_2) \int_s^t d\Gamma_{a_1}(x_3)\\ &&-\int_s^t d\Gamma_{a_2}(x_1) \int_s^{x_1}
d\Gamma_{a_3}(x_2) \int_s^{x_1} d\Gamma_{a_1}(x_3)\\ &=&I^{ts}_\Gamma(\tddeux{$a_2$}{$a_3$}\:
\tdun{$a_1$}\:)-I^{ts}_\Gamma(\hspace{2mm}\tdtroisun{$a_2$}{$a_1$}{\hspace{-1.5mm}$a_3$}\:),
\end{eqnarray*} so $\T^{\sigma}=\tddeux{$a_2$}{$a_3$}\:
\tdun{$a_1$}\:-\hspace{2mm}\tdtroisun{$a_2$}{$a_1$}{\hspace{-1.5mm}$a_3$}$. Note that all permutation
graphs $\T^{\sigma}$ with $\sigma$ fixed are obtained from  the same sum of heap-ordered forests (also
denoted by $\T^{\sigma}$, by abuse of notation) by including
 the decorations of $\mathfrak{T}_n$ permuted
by $\sigma$.

\medskip

As an elementary Corollary of Definition \ref{def:splitting} and  Proposition
\ref{prop:permutation-graph}, one obtains:

\begin{Corollary}[Fourier normal ordering for smooth paths]

Let $\Gamma$ be a smooth path and $\mathfrak{T}_n\in {\bf Sh}$ a trunk tree with $n$ vertices and
decoration $\ell$, then \BEQ I^{ts}_{\Gamma}(\mathfrak{T}_n)=\sum_{\sigma\in \Sigma_n}
I^{ts}_{\mu^{\sigma}}(\T^{\sigma}).\EEQ

\end{Corollary}


\subsection{Fourier normal ordering and regularization}


Formal rough paths over $\Gamma$ will be reconstructed out of {\em tree data} $\phi^t_{\T}$
defined arbitrarily for each tree $\T$, and then extended by multiplication to forests, as we shall now see.

\begin{Definition} \label{def:barchi-char}

\begin{itemize}

\item[(i)] For every heap-ordered $\T$ with $n$ vertices, and $t\in\R$,
 let $\phi^t_{\T}: {\cal P}^{\T} Meas(\R^n)\to\R, \mu\mapsto \phi^t_{\mathbb{T}}(\mu)$, also written
 $\phi^t_{\mu}(\mathbb{T})$
 be a family of linear forms such that:

\begin{itemize} \item[(a)]
$\phi^t_{d\Gamma(i)}(\T_1)-\phi^s_{d\Gamma(i)}(\tdun{i})=I^{ts}_{\Gamma}(\T_1)=\Gamma_t(i)-\Gamma_s(i)$
if $\T_1$ is the trivial heap-ordered tree with one vertex;

\item[(b)]  if $\mathbb{T}_i$, $i=1,2$ are heap-ordered trees with $n_i$ vertices, and
    $\mu_i\in
 {\cal P}^{\mathbb{T}_i} Meas(\R^{n_i})$, $i=1,2$, the following
multiplicative property holds, \BEQ \phi^t_{\mu_1}(\mathbb{T}_1)
\phi^t_{\mu_2}(\mathbb{T}_2)=\phi^t_{\mu_1\otimes\mu_2}(\mathbb{T}_1 \wedge\mathbb{T}_2),
\label{eq:T1wedgeT2}\EEQ where $\T_1\wedge\T_2$ is the non-decorated product $\T_1.\T_2$ with
labels of $\T_2$ shifted by $n_1$ \footnote{The
 product $\T_1\wedge\T_2$ defines actually the product of the Hopf algebra of heap-ordered trees
 \cite{FoiUnt}};

\item[(c)]  ({\em naturality property})
 the following invariance condition under reindexation of the vertices
holds, see preceding two Remarks, \BEQ \phi^t_\mu(\T)=\phi^t_{\mu \circ
\sigma}(\sigma^{-1}.\mathbb{F}) \EEQ if  $\sigma$ -- which acts by permuting  the vertices of
$\T$ -- is such that \BEQ (i\twoheadrightarrow j\ {\mathrm{in}}\ \T)\Rightarrow
(i\twoheadrightarrow j\ {\mathrm{in}}\ \sigma^{-1}(\T)).\EEQ

\end{itemize}

\item[(ii)] Let, for $\Gamma=(\Gamma(1),\ldots,\Gamma(d))$,
 ${\chi}_{\Gamma}^t:{\bf Sh}\to\R$ be the linear form on ${\bf Sh}$ defined by
\BEQ {\chi}_{\Gamma}^t (\mathfrak{T}_n):=\sum_{\sigma\in \Sigma_n}
\phi^t_{\mu_{(\Gamma,\ell)}^{\sigma}}(\mathbb{T}^{\sigma}), \quad
\mathfrak{T}_n=(\ell(1)\ldots\ell(n))\EEQ as in Proposition \ref{prop:permutation-graph}.
\end{itemize} \end{Definition}

\bigskip

The main result is the following.

\begin{Proposition}[rough path construction by Fourier normal ordering]  \label{lem:barchi-char}
 For every path $\Gamma$ such that $\chi_{\Gamma}^t$ is well-defined,
 ${\chi}_{\Gamma}^t$ is a character of ${\bf Sh}$.

Consequently, the following formula  for $\mathfrak{T}_n\in {\bf Sh}$,  $n\ge 1$, with $n$ vertices and
decoration $\ell$,
  \BEQ J^{ts}_{\Gamma}(\ell(1),\ldots,\ell(n)):=
   {\chi}^t_{\Gamma} \ast ({\chi}^s_{\Gamma}\circ S)
(\mathfrak{T}_n) \label{eq:Jchi}  \EEQ defines a {\em formal} rough path over $\Gamma$.

Furthermore, the following equivalent definition holds, \BEQ
J^{ts}_{\Gamma}(\mathfrak{T}_n):=\sum_{\sigma\in\Sigma_n} \left( \phi^t \ast (\phi^s\circ
\bar{S})\right)_{\mu^{\sigma}_{(\Gamma,\ell)}}(\mathbb{T}^{\sigma}), \label{eq:Jphi}  \EEQ where the
convolution in the right equation is defined by reference to the (heap-ordered) tree coproduct, namely,
one sets \BEQ (\phi^t\ast(\phi^s\circ\bar{S}))_{\nu}(\T)=\sum_{\vec{v}\models V(\T)} \phi^t_{
\otimes_{v\in V(\Roo_{\vec{v}}\T)} \nu_v} (\Roo_{\vec{v}}\T)  \phi^s_{ \otimes_{v\in
V(\Lea_{\vec{v}}\T)} \nu_v} (\bar{S}(\Lea_{\vec{v}}\T)) \EEQ for a tensor measure
$\nu=\nu_1\otimes\ldots\otimes\nu_n$, and by multilinear extension

\BEA &&  \left(\phi^t\ast(\phi^s\circ\bar{S})\right)_{\nu}(\T)=(2\pi)^{-n/2} \int {\cal
F}\nu(\xi_1,\ldots,\xi_n)d\xi_1\ldots d\xi_n \ \cdot\nonumber\\ && \cdot \sum_{\vec{v}\models V(\T)}
\phi^t_{\otimes_{v\in V(Roo_{\vec{v}}(\T))} e^{\II  x_v\xi_v} dx_v}(\Roo_{\vec{v}}\T)
\phi^s_{\otimes_{v\in V(Lea_{\vec{v}}(\T))} e^{\II  x_v\xi_v} dx_v}(\bar{S}(\Lea_{\vec{v}}\T)) ,\quad
\T\in {\cal F}_{ho}(n).  \nonumber\\ \label{convol} \EEA

for an arbitrary measure $\nu\in Meas(\R^n)$, $\nu=(2\pi)^{-n/2} \int d\vec{\xi} {\cal F}\nu(\vec{\xi})
\otimes_{j=1}^n e^{\II x_j\xi_j} dx_j.$

\end{Proposition}

Quite naturally, we shall call  formula (\ref{eq:Jchi}), resp. (\ref{eq:Jphi}) the {\em shuffle
convolution}, resp. {\em tree convolution definition} of $J$. \medskip

Assuming $\Gamma$ is smooth, then defining $\phi^t$ as the skeleton integral $\SkI^t$ yields trivially
by recombination $\chi^t_{\Gamma}=\SkI^t$ too, and then $J_{\Gamma}^{ts}=I^{ts}_{\Gamma}$ is the
canonical rough path over $\Gamma$.
 Proposition \ref{lem:barchi-char} shows
that the same recombination algorithm yields a rough path over $\Gamma$ whenever $\phi^t$ satisfies
conditions (a), (b) and (c) of Definition \ref{def:barchi-char}. It is actually clear from Definition
\ref{def:barchi-char} that {\em any} rough path over $\Gamma$ may be  obtained in this way
\cite{FoiUnt}.

\bigskip

The enormous advantage now with respect to the original problem is that one may construct as many linear
forms $\phi^t$ as one wishes by assigning some arbitrary value to $\phi^t_{\mu}(\T)$, $\T$ ranging over
all (heap-ordered) trees with $\ge 2$ vertices, and extending to forests by multiplication following
condition (b).  

\medskip

It is now natural to try and define $\phi^t$ as some regularized skeleton integral in such a way that
$J^{ts}_{\Gamma}$ satisfies the H\"older continuity property (i) in Definition \ref{def:rough-path}. We
shall do so in the next sections by {\em renormalizing} skeleton integrals.

\bigskip

For the sequel,  we shall start from the {\em tree convolution definition} (\ref{eq:Jphi}) of $J$, which
will be used in the following guise.  Assume $\nu=\mu^{\sigma}_{(\Gamma,\ell)}$ and $\T=\T_1\wedge\ldots
\wedge\T_p$ is the (heap-ordered) product of $p$ trees. Set $\hat{\nu}_{\T'}(\vec{\xi})=\otimes_{v\in
V(\T')} {\cal F}(\Gamma'(\ell\circ\sigma(v)))(\xi_v) e^{\II x_v \xi_v} dx_v$ for $\T'$ subtree of $\T$
and $\vec{\xi}=(\xi_v)_{ v\in V(\T')}$.  Then, by the multiplicative property (b) for $\phi^s$ and
$\phi^t$, see Definition \ref{def:barchi-char}, \BEQ (\phi^t\ast (\phi^s\circ\bar{S}))_{\nu}(\T) =
(2\pi)^{-n/2} \int d\xi_1\ldots d\xi_n    {\bf 1}_{|\xi_1|\le\ldots\le|\xi_n|} \prod_{q=1}^p (\phi^t\ast
(\phi^s\circ\bar{S}))_{\hat{\nu}_{\T_q}((\xi_v)_{v\in V(\T_q)})}(\T_q). \label{eq:1.35} \EEQ Now the
inductive definition of the antipode implies \BEA
&&(\phi^t\ast(\phi^s\circ\bar{S}))_{\hat{\nu}_{\T_q}((\xi_v)_{v\in V(\T_q)})}(\T_q) \nonumber\\
&&=\phi^t_{\hat{\nu}_{\T_q}((\xi_v)_{v\in V(\T_q)})}(\T_q)+\phi^s_{\hat{\nu}_{\T_q}((\xi_v)_{v\in
V(\T_q)})}(\bar{S}(\T_q)) \nonumber\\ && \qquad  +\sum_{\vec{v}\in V(\T_q)}
\phi^t_{\hat{\nu}_{\Roo_{\vec{v}}\T_q}((\xi_v)_{v\in V(\Roo_{\vec{v}}\T_q)})}(\Roo_{\vec{v}}\T_q)
\phi^s_{\hat{\nu}_{\Lea_{\vec{v}}\T_q}((\xi_v)_{v\in
V(\Lea_{\vec{v}}\T_q)})}(\bar{S}(\Lea_{\vec{v}}\T_q)) \nonumber\\
&&=(\phi^t-\phi^s)_{\hat{\nu}_{\T_q}((\xi_v)_{v\in V(\T_q)})}(\T_q) \nonumber\\ && \quad
+\sum_{\vec{v}\models V(\T_q),\vec{v}\not=\emptyset}
 (\phi^t-\phi^s)_{\hat{\nu}_{\Roo_{\vec{v}}\T_q}((\xi_v)_{v\in
 V(\Roo_{\vec{v}}\T_q)})}(\Roo_{\vec{v}}\T_q)
\phi^s_{\hat{\nu}_{\Lea_{\vec{v}}\T_q}((\xi_v)_{v\in
V(\Lea_{\vec{v}}\T_q)})}(\bar{S}(\Lea_{\vec{v}}\T_q)) \nonumber\\ \EEA

Finally, applying iteratively the inductive definition of the antipode leads to an expression of
$\bar{S}(\Lea_{\vec{v}}\T_q)$ in terms of a sum of forests obtained by multiple cuts an in
\cite{ConKre00}. Applying once again the multiplicative property to $\phi^s$ yields
$(\phi^t\ast(\phi^s\circ \bar{S}))_{\hat{\nu}_{\T_q}((\xi_v)_{v\in V(\T_q)})}(\T_q)$ as a sum of terms
of the form \BEA &&  \Phi^{ts}(\T_q;\vec{\xi}; \vec{v},(\T'_j)):= \nonumber\\ &&\qquad
(\phi^t-\phi^s)_{\hat{\nu}_{\Roo_{\vec{v}}\T_q}((\xi_v)_{v\in
V(\Roo_{\vec{v}}\T_q)})}(\Roo_{\vec{v}}\T_q) \prod_{j=1}^J \phi^s_{\hat{\nu}_{\T'_j}((\xi_v)_{v\in
V(\T'_j)})}(\T'_j), \nonumber \label{eq:1.37}\EEA
with $V(\T_q)=V(\Roo_{\vec{v}}\T_q)\cup \uplus_{j=1}^J
V(\T'_j).$


\section{Feynman diagram reformulation}


Let $\T$ be a forest. We shall show in this section how to  compute tree  skeleton  integrals
$\SkI_B(\T)$ of fractional Brownian motion by means of Feynman diagrams of a particular type.
Computations are based on the {\em harmonizable representation} of fBm,
\BEQ B_t(i)=(2\pi c_{\alpha})^{-\half} \int \frac{e^{\II t\xi}-1}{\II\xi} |\xi|^{\half-\alpha}
dW_{\xi}(i), \quad 1\le i\le d \label{eq:harm} \EEQ
where $(W_{\xi}(1),\ldots,W_{\xi}(d))$ are $d$ independent, identically distributed complex Brownian motions such that $W_{-\xi}(i)=\overline{W_{\xi}(i)}$. With the usual normalization choice
$(2\pi c_{\alpha})^{-\half}=\half\sqrt{-\frac{\alpha}{\cos\pi\alpha\Gamma(-2\alpha)}}$, $(B_t)_{t\in\R}$
is the unique centered Gaussian process with covariance
\BEQ \esper B_s B_t=\half(|s|^{2\alpha}+|t|^{2\alpha}-|t-s|^{2\alpha}).\EEQ

 Quite
generally, the associated physical theory contains particles of 2 types, corresponding to two Gaussian
fields, $\sigma$, resp. $\phi$, whose propagators are represented by simple, resp. double lines.
Vertices are of type $(\phi\sigma^n)_{n\ge 2}$, namely, at each vertex meet $n\ge 2$
 simple lines and
exactly 1 double line. More specifically, we shall only need to consider {\em tree Feynman diagrams} in
an unusual sense, namely, Feynman diagrams such that the subset of {\em simple} lines contains no
loops.

We shall also speak for convenience of Feynman {\em half-diagrams}, which are Feynman diagrams in the
above sense, except that it also possibly admits -- besides {\em true} external $\phi$-legs -- {\em
uncontracted  $\phi$-legs}, which  are assumed to be cut in the middle (this implies special evaluation
rules as we shall see). On the other hand, {\em contracted $\phi$-legs} are always internal lines. Gluing a Feynman half-diagram $G^{\half}$ with its image in a mirror along the
middle of its external double lines yields a {\em symmetric Feynman diagram} $G=(G^{\half})^2$.

\begin{figure}[h]
  \centering
   \includegraphics[scale=0.35]{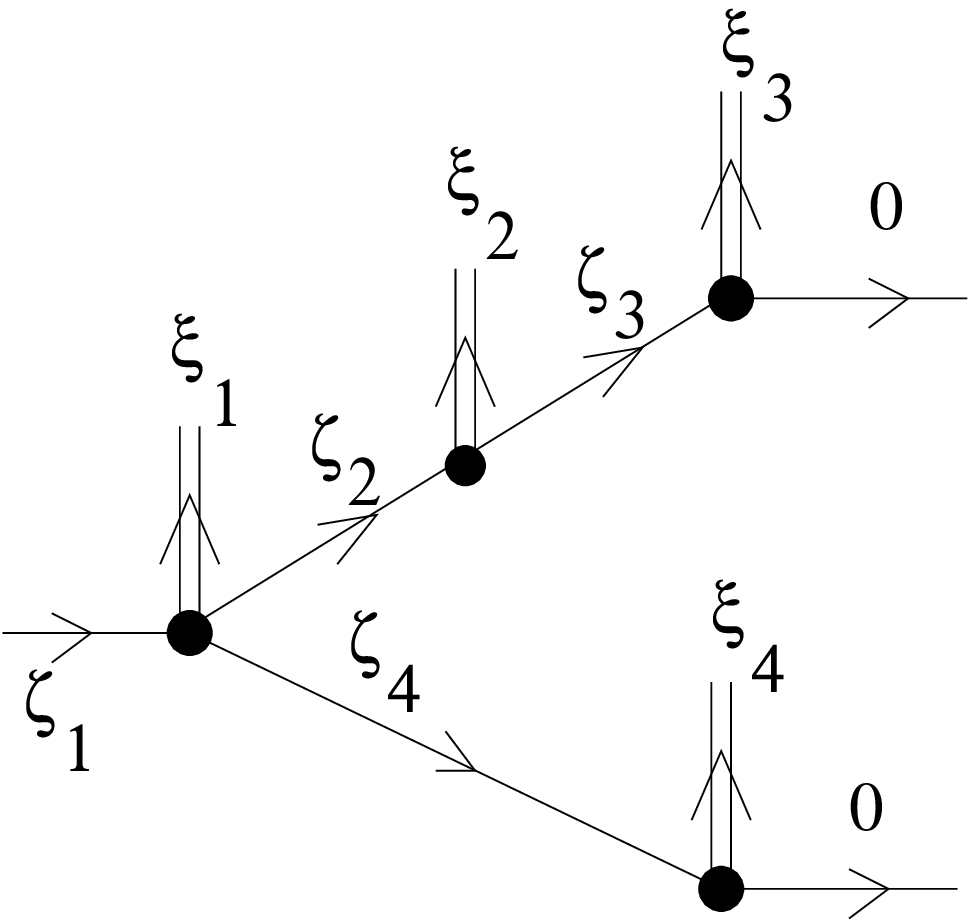}
   \caption{\small{Feynman half-diagram $G^{\half}(\T)$ associated to $\T$}}
  \label{Fig-renbis1}
\end{figure}

\begin{figure}[h]
  \centering
   \includegraphics[scale=0.45]{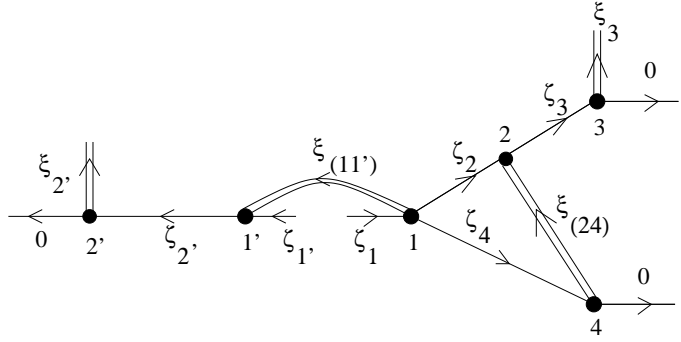}
   \caption{\small{Feynman half-diagram $G^{\half}(\T'.\T;(11'),(24))$
 associated to $\T'.\T$. By momentum conservation,
 $\xi_{2'}+\xi_3=\zeta_{1'}+
\zeta_1.$}}
  \label{Fig-ren1contract}
\end{figure}

\begin{figure}[h]
  \centering
   \includegraphics[scale=0.35]{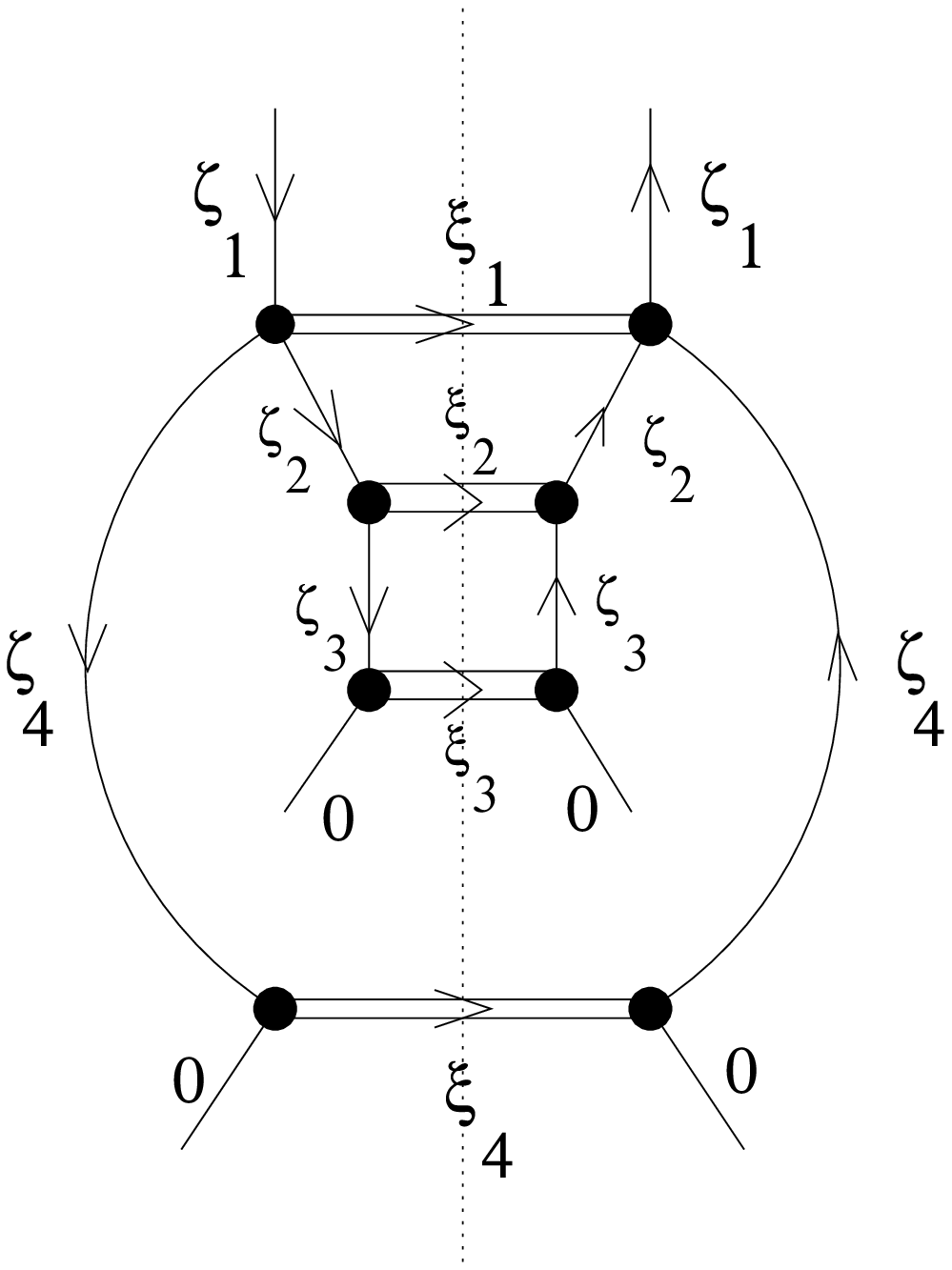}
   \caption{\small{Full Feynman diagram $G(\T)=(G^{\half}(\T))^2$ associated to  Fig.
   \ref{Fig-renbis1}.}}
  \label{Fig-renbis2}
\end{figure}

\begin{figure}[h]
  \centering
   \includegraphics[scale=0.55]{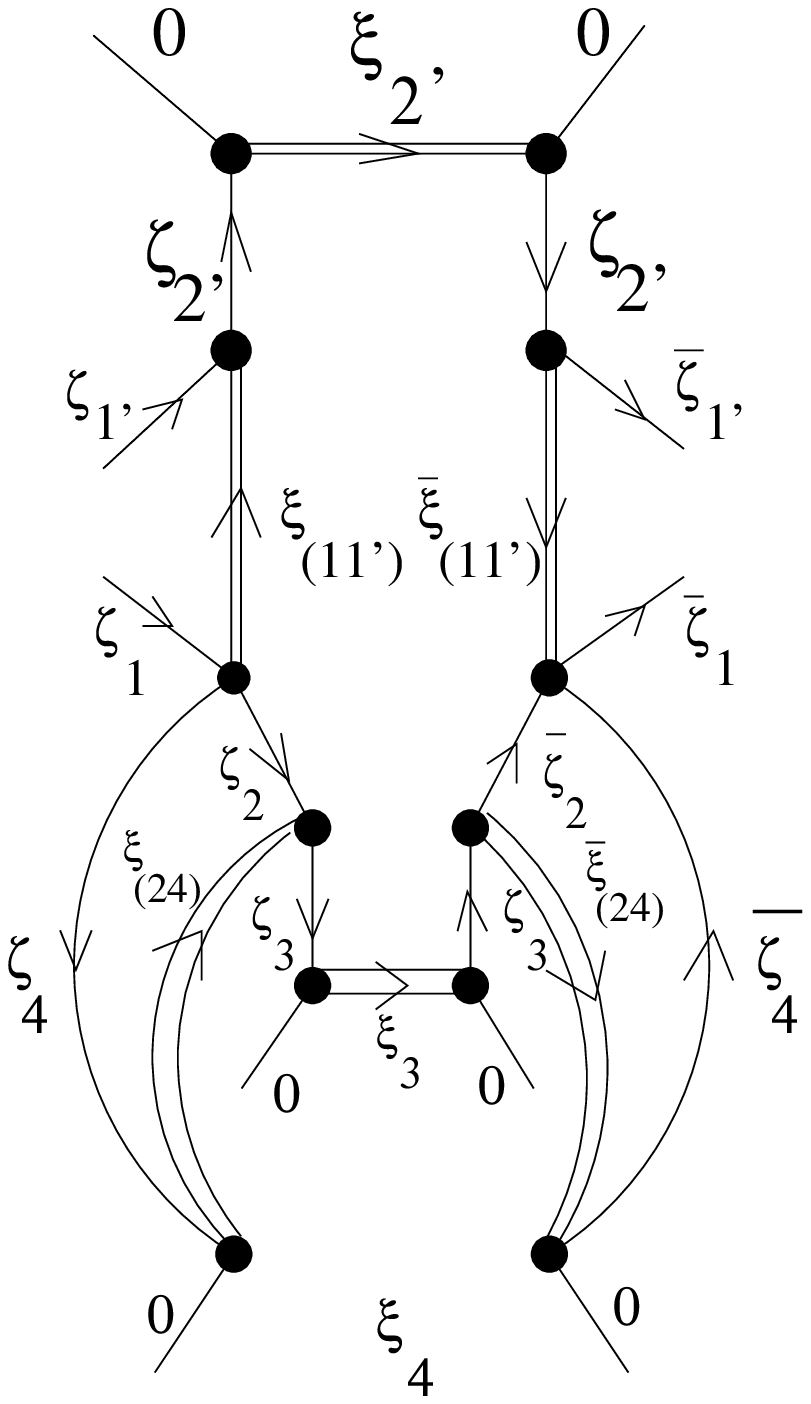}
   \caption{\small{Full Feynman diagram $G(\T'.\T;(1 1'),(24))=(G^{\half}(\T'.\T;(11'),(24)))^2$
associated to Fig. \ref{Fig-ren1contract}.}}
  \label{Fig-renbis2contract}
\end{figure}

A tree (or more generally a forest) $\T$ determines a unique tree Feynman half-diagram $G^{\half}(\T)$
(called: {\em uncontracted tree Feynman half-diagram associated to $\T$}), admitting only uncontracted
$\phi$-legs, whose underlying tree
structure of simple lines is that of $\T$, see Fig. \ref{Fig-renbis1} and
Fig. \ref{Fig-renbis2}.  One always assigns {\em zero momentum} to the simple
external lines attached to the {\em leaves} of $\T$. All other tree Feynman half-diagrams are obtained
from some $G^{\half}(\T)$ by pairwise contracting some of the uncontracted $\phi$-legs, and denoted
accordingly, see
  Fig. \ref{Fig-ren1contract} and
Fig. \ref{Fig-renbis2contract}.  In these diagrams, $\T=\tdquatredeux{1}{4}{2}{3}$ and
$\T'=\tddeux{1'}{2'}$.  If $G^{\half}$ is a tree Feynman half-diagram (in the same sense as for Feynman
diagrams) then G is called a {\em symmetric tree Feynman diagram}. By definition, a symmetric tree
Feynman diagram has no external double line.

Let us now  define  Feynman rules. If $G$ is a diagram or half-diagram, the set of vertices, resp.
internal lines shall be denoted by $V(G)$, resp. $L(G)$. The set of external lines is
denoted by $L_{ext}(G)$.

\begin{Definition}[Feynman rules]

Let $G$ be a Feynman diagram or half-diagram consisting of simple lines, double lines and vertices $v$
connecting one double line with $2,3,\ldots$ double lines, and a certain number of external
lines.
 Each line is oriented and decorated by
a real-valued momentum, conventionally denoted by $\zeta_i$, resp. $\xi_i$  or $\xi_{(ij)}$ for some
index $i$ or pair contraction $(ij)$  for simple, resp. double lines; reversing the orientation of the
line is equivalent to changing the sign of the momentum. The momentum preservation relation holds,
namely, the sum of all momenta at any vertex is zero.
 We denote by $I_{\sigma}(G)$, resp. $I_{\phi}(G)$, the number of internal simple, resp. double lines,
 so
that simple, resp. double lines may be thought as propagators of some field denoted by $\sigma$, resp.
$\phi$. We also let $I(G):=I_{\sigma}(G)+I_{\phi}(G)$ be the total number of internal lines.

\begin{itemize} \item[(i)] Feynman half-diagrams

Let $G^{\half}$ be a half-diagram. Associate $\zeta_i^{-1}$ to each internal simple line with
momentum $\zeta_i$,
 $|\xi_i|^{\half-\alpha}$ to each uncontracted $\phi$-leg with momentum $\xi_i$, and
 $|\xi_{(ij)}|^{1-2\alpha}$ to each
contracted double line with momentum $\xi_{(ij)}$. The result is a function of the momenta of the
external lines, $\mathbf{\zeta}_{ext}$ and $\mathbf{\xi}_{ext}$, denoted by
$A_{G^{\half}}(\mathbf{\zeta}_{ext},\mathbf{\xi}_{ext}).$ We shall denote by  $\zeta_{ext}$, resp.
$\xi_{ext}$ the sum of the momenta of the external simple, resp. double lines.

In the particular case when $G^{\half}=G^{\half}(\T;(i_1 i_2),\ldots,(i_{2p-1} i_{2p}))$, $p\ge 0$,
comes from a tree, one has $\mathbf{\xi}_{ext}=\{\xi_v\ |\ \xi_v \ {\mathrm{uncontracted}} \}$, and
\BEA &&  A_{G^{\half}}(\mathbf{\zeta}_{ext},\mathbf{\xi}_{ext})= \del(\zeta_{ext}+\xi_{ext})\nonumber\\ &&
\quad \prod_{v\in V(\T)\ |\ \xi_v \ {\mathrm{uncontracted}}}
 |\xi_v|^{\half-\alpha}. \int \prod_{v\in V(\T)\setminus\{ {\mathrm{roots}} \}} \frac{1}{\zeta_v}\
 \cdot\
\prod_{q=1}^p |\xi_{(i_{2q-1} i_{2q})}|^{1-2\alpha} d\xi_{(i_{2q-1} i_{2q})}, \nonumber\\
\label{eq:2.1}  \EEA all external $\zeta$-momenta attached to the leaves of $\T$ vanishing, as
explained before.

\item[(ii)] Feynman diagrams

Associate $\zeta^{-1}_i$ to each internal simple line with momentum $\zeta_i$,
 and $|\xi_i|^{1-2\alpha}$, resp. $|\xi_{(ij)}|^{1-2\alpha}$  to each internal double line with
 momentum $\xi_i$, resp. $\xi_{(ij)}$.

The resulting amplitude of the amputated diagram (i.e., shorn of its external legs),
 function of the momenta of the external lines, $\mathbf{\zeta}_{ext}$ and
$\mathbf{\xi}_{ext}$,  is denoted by  $A_{G}(\mathbf{\zeta}_{ext},\mathbf{\xi}_{ext})$.

In the particular case when $G=(G^{\half})^2$ is a  symmetric tree Feynman diagram, denoting by
$\bar{\zeta}_i$, $\bar{\xi}_{(ij)}$ the momenta of the mirror lines,  and by $\zeta_{ext}$, resp. $\bar{\zeta}_{ext}$ the sum of the external $\zeta$-, resp.
$\bar{\zeta}$-momenta, one has \BEA  &&  A_G(\mathbf{\zeta}_{ext};\bar{\mathbf{\zeta}}_{ext})=
\nonumber\\ && \quad  \int \prod_{v\in V(\T) \ |\ \xi_v \ {\mathrm{uncontracted}}} d\xi_v
 \left(A_{G^{\half}}(\mathbf{\zeta}_{ext},\mathbf{\xi}_{ext})\right)
 \left(A_{G^{\half}}(\bar{\mathbf{\zeta}}_{ext},\mathbf{\xi}_{ext})\right),
\nonumber\\ \label{eq:2.2} \EEA where $\mathbf{\xi}_{ext}=\{\xi_v;\ \xi_v\ {\mathrm{uncontracted}}\}$
as in (i).  Eq. (\ref{eq:2.2}) contains a hidden $\del$-function $\del(\zeta_{ext}+\bar{\zeta}_{ext})$ due to overall momentum conservation.
\end{itemize} \end{Definition}

{\bf Remark.} The relation between the $\zeta$- and $\xi$-coordinates for a half- or full diagram coming
from a tree is simply $\zeta_v=\xi_v+\sum_{w\twoheadrightarrow v} \xi_w$ or conversely,
$\xi_v=\zeta_v-\sum_{w\to v} \zeta_w$, where $\{w:w\twoheadrightarrow v\}$, resp. $\{w:w\to v\}$ are the descendants, resp. children of $v$, and one has set
$\xi_i=-\xi_j=\xi_{(ij)}$ for contracted double lines.

\bigskip

Let $G$ be a connected Feynman diagram. It contains $I(G)-|V(G)|+1$ independent momenta: namely, there is
one momentum constraint at each vertex, which gives altogether $|V(G)|-1$ independent constraints,
because  the global translation invariance has already been taken into account by demanding that the sum
of the external momenta be zero. Remove one internal line at each vertex, so that all remaining momenta
are independent. The set of all lines which have been removed, together with the vertices at the end of
the lines, constitute a subdiagram of $G$ with no loops, hence a sub-forest. For such a choice of lines,
$L'(G)$, say,  we let  $(z_{\ell})_{\ell\in L(G)\setminus L'(G)}$, $z=\zeta$ or $\xi$, be the set of
remaining, independent momenta. Each $z_{\ell'}$, $\ell'\in L'(G)$, may be written uniquely as some
linear combination
 $z_{\ell'}=z_{\ell'}\left( (z_{\ell})_{\ell \in L_{ext}(G)\cup (L(G)\setminus L'(G))} \right)$,
 which yields an explicit formula for $A_G$,
\BEQ A_G(\mathbf{z}_{ext})=\del(z_{ext})\int \prod_{\ell\in L(G)\setminus L'(G)} dz_{\ell} \ .\
\prod_{\ell\in L_{\phi}(G)} |\xi_{\ell}|^{1-2\alpha} \prod_{\ell\in L_{\sigma}(G)} \zeta_{\ell}^{-1},
\label{eq:2.3} \EEQ where $z_{ext}$ is the sum of the external momenta.

The relation with iterated integrals of fractional Brownian motion is the following.

\begin{Lemma} \label{lem:int-Feynman}

\begin{enumerate} \item Let $\T$ be a tree with $n$ vertices and root indexed by $1$. Then, see
 \ref{eq:1.16}and \ref{eq:harm},
 \BEA  &&
\overline{\SkI}^t_{B}(\T)=  (i\sqrt{2\pi c_{\alpha}})^{-n} \int \prod_{v\in V(\T)} dW_{\xi_v}(\ell(v))
\nonumber\\ && \qquad \int \frac{e^{\II t\zeta_1}}{\zeta_1}  d\zeta_1
A_{G^{\half}(\T)}(\mathbf{\zeta}_{ext}=(\zeta_1,0),\mathbf{\xi}_{ext}=(\xi_v)_{v\in V(\T)} ), \EEA with \BEQ
A_{G^{\half}(\T)}(\mathbf{\zeta}_{ext}=(\zeta_1,0),\mathbf{\xi}_{ext}= (\xi_v)_{v\in V(\T)}
)=\del(\zeta_1+\xi_{ext}) \prod_{v\in V(\T)} |\xi_v|^{\half-\alpha} \ \cdot\ \prod_{v\in
V(\T)\setminus\{1\}} \frac{1}{\zeta_v}.\EEQ Hence, {\em provided} the decorations $(\ell(v))_{v\in
V(\T)}$ are all distinct, \BEQ   \Var (\overline{\SkI}^t_B(\T) - \overline{\SkI}^s_B(\T) )=(2\pi
c_{\alpha})^{-n} \int \frac{d\zeta_1}{\zeta_1^2} |e^{\II t\zeta_1}-e^{\II s\zeta_1}|^2
A_{G(\T)}(\mathbf{\zeta}_{ext} =(\zeta_1,0)), \EEQ with \BEA  A_{G(\T)}(\mathbf{\zeta}_{ext}=(\zeta_1,0)) &=&
 \int \prod_{v\in V(\T)} d\xi_v \left| A_{G^{\half}}((\zeta_1,0),\mathbf{\xi})\right|^2
 \label{eq:2.6,5} \\
 &=&\int \prod_{v\in V(\T)} d\xi_v
\del(\zeta_1+\sum_{v\in V(\T)} \xi_v) \prod_{v\in V(\T)} |\xi_v|^{1-2\alpha} \ \cdot\ \prod_{v\in
V(\T)\setminus\{1\}} \frac{1}{\zeta_v^2}. \nonumber\\ \EEA

\item Let more generally $\T=\T_1\ldots \T_q$ and $\T'=\T'_1\ldots\T'_{q'}$, $q,q'\ge 1$, with roots
    $r_1,\ldots,r_q$,$r'_1,\ldots,r'_{q'}$. Consider some multiple contraction $(i_1
    i_2),\ldots,(i_{2p-1} i_{2p})$ of $\prod_{m=1}^q
    (\overline{\SkI}^t_B(\T_{m})-\overline{\SkI}^s_B(\T_{m})) \prod_{m'=1}^{q'}
    \overline{\SkI}^s_B(\T'_{m'})$ connecting the vertices of $\T$ and $\T'$,
 which we write for short \\ $\del \overline{\SkI}^{ts}_B \overline{\SkI}^s_B (\T,\T';
(i_1 i_2),\ldots,(i_{2p-1}i_{2p}))$, and let $G^{\half}:= G^{\half}(\T.\T';(i_1
i_2),\ldots,(i_{2p-1} i_{2p}))$ be the corresponding (connected) Feynman half-diagram. Then

 \BEA &&  \del\overline{\SkI}^{ts}_{B}\overline{\SkI}^s(\T,\T';(i_1 i_2),\ldots,(i_{2p-1} i_{2p}))=
 (\II\sqrt{2\pi c_{\alpha}})^{-(|V(\T)|+|V(\T')|)}
\nonumber\\ && \quad \int \prod_{v\in V(\T.\T')\ |\ \xi_v\ {\mathrm{uncontracted}}}
dW_{\xi_v}(\ell(v)) \prod_{m=1}^q \frac{e^{\II t\zeta_{r_m}}-{e^{\II s\zeta_{r_m}}
}}{\zeta_{r_m}} d\zeta_{r_m} \nonumber\\ && \qquad
 \prod_{m'=1}^{q'} \frac{e^{\II s\zeta_{r'_{m'}}}}{\zeta_{r'_{m'}}} d\zeta_{r'_{m'}}
 A_{G^{\half}}
(\mathbf{\zeta}_{ext}=((\zeta_{r_m}),(\zeta_{r'_{m'}}),0),\mathbf{\xi}_{ext}) \nonumber\\ \EEA where
$\mathbf{\xi}_{ext}:=\{(\xi_v)_{v\in V(\T)}\ |\ \xi_v \ {\mathrm{uncontracted}}\}.$

Assume furthermore all non-contracted indices $\ell(i)$, $i\not=i_1,\ldots,i_{2p}$ are distinct.
Then
 \BEA &&  \Var  \del\overline{\SkI}^{ts}_{B}\overline{\SkI}^s(\T,\T';(i_1 i_2),\ldots,(i_{2p-1}
 i_{2p}))= (2\pi c_{\alpha})^{-(|V(\T)|+|V(\T')|)} \nonumber\\
&& \quad \int \prod_{m=1}^q d\zeta_{r_m} d\bar{\zeta}_{r_m} \left( \frac{e^{\II t\zeta_{r_m}}-e^{\II
s\zeta_{r_m}}}{\zeta_{r_m}}\right)  \left( \frac{e^{\II t\bar{\zeta}_{r_m}}-e^{\II
s\bar{\zeta}_{r_m}}}{\bar{\zeta}_{r_m}}\right)   \nonumber\\
 && \qquad \int \prod_{m'=1}^{q'} \frac{d\zeta_{r'_{m'}} d\bar{\zeta}_{r'_{m'}}}{ \zeta_{r'_{m'}}
 \bar{\zeta}_{r'_{m'}}}
A_{(G^{\half})^2} (\mathbf{\zeta}_{ext}=((\zeta_{r_m}),(\zeta_{r'_{m'}}),0);
\bar{\mathbf{\zeta}}_{ext}=((\bar{\zeta}_{r_m}),(\bar{\zeta}_{r'_{m'}}),0) ). \nonumber\\ \EEA
\end{enumerate} \end{Lemma}

{\bf Example.} Let $\T$, $\T'$  be as in Fig.  \ref{Fig-ren1contract}, \ref{Fig-renbis2contract}. Then:

\BEA && \del\overline{\SkI}^{ts}_{B}\overline{\SkI}_B^s(\T,\T';(1 1'),(2 4)) \nonumber\\ && =
(\II\sqrt{2\pi
c_{\alpha}})^{-6} \int dW_{\xi_{2'}}(\ell(2'))dW_{\xi_3}(\ell(3)) \left[\frac{e^{\II t\zeta_1}-e^{\II
s\zeta_1}}{\zeta_1} d\zeta_1\right] \left[\frac{e^{\II s\zeta_{1'}}}{\zeta_{1'}}
d\zeta_{1'}\right] A_{G^{\half}}(\zeta_1,\zeta_{1'},\xi_{2'},\xi_3)  \nonumber\\ \EEA with \BEQ
A_{G^{\half}}(\zeta_1,\zeta_{1'},\xi_{2'},\xi_3)=\del(\zeta_1+\zeta_{1'}+\xi_{2'}+\xi_3)
|\xi_{2'}\xi_3|^{\half-\alpha} \frac{|\xi_{(11')}\xi_{(24)}|^{1-2\alpha} d\xi_{(11')}
d\xi_{(24)}}{(\xi_{(24)}+ \xi_3)\xi_3 \xi_{(24)}\xi_{2'}}.\EEQ

As for its variance, assuming $\ell(2')\not=\ell(3)$,

\BEA && \Var \del\overline{\SkI}^{ts}_{B}\overline{\SkI}^s(\T,\T';(1 1'),(2 4)) \nonumber\\ && =(2\pi
c_{\alpha})^{-6} \int \frac{d\zeta_1 d\zeta_{1'} d\bar{\zeta}_1 d\bar{\zeta}_{1'}}{\zeta_{1'}
\bar{\zeta}_{1'}}  \left[\frac{e^{\II t\zeta_1}-e^{\II s\zeta_1}}{\zeta_1}\right]
\left[\frac{e^{\II t\bar{\zeta}_1}-e^{\II s\bar{\zeta}_1}}{\bar{\zeta}_1}\right]
 A_{(G^{\half})^2}(\zeta_1,\zeta_{1'};\bar{\zeta}_1,\bar{\zeta}_{1'})
 \nonumber\\ \EEA
with \BEA &&  A_{(G^{\half})^2}(\zeta_1,\zeta_{1'};\bar{\zeta}_1,\bar{\zeta}_{1'}) = \int
d\xi_{2'}d\xi_3 d\xi_{(24)} d\xi_{(11')} d\bar{\xi}_{(24)} d\bar{\xi}_{(11')}
\del(\zeta_1+\zeta_{1'}+\xi_{2'}+\xi_3) \nonumber\\ && \quad \del(\bar{\zeta}_1+
\bar{\zeta}_{1'}+\xi_2+\xi_3)
  |\xi_{2'}\xi_3|^{1-2\alpha} \frac{ |\xi_{(24)}\xi_{(11')}
  \bar{\xi}_{(24)}\bar{\xi}_{(11')}|^{1-2\alpha}}{ (\xi_3\xi_{2'})^2
  (\xi_{(24)}+
\xi_3) (\bar{\xi}_{(24)}+\xi_3)    \xi_{(24)} \bar{\xi}_{(24)} }. \nonumber\\ \EEA


\section{Definition of renormalization scheme}


We present here the general features of the BPHZ renormalization scheme, together with its multi-scale
formulation which will allow us to prove H\"older regularity. It relies

\begin{itemize} \item[(i)] on the choice of a set of graphs called {\em diverging graphs}. In general
(see subsection 3.1 below) it is simply the subset of Feynman graphs $G$ such that $\omega(G)>0$, where $\omega$ is the
{\em overall degree of divergence} (or simply degree of homogeneity) of the graph.

\item[(ii)] on a choice of {\em regularization scheme}. Here we choose the Taylor evaluation at zero
    external momenta, denoted by $\tau$. To be definite, if $A_g(z_{ext,1},\ldots, z_{ext,N_{ext}})$
    is the amplitude of the graph $g$ with $N_{ext}$ external momenta, then $\tau_g
    A_g(z_{ext,1},\ldots,z_{ext,N_{ext}})=A_g(0,\ldots,0).$

\end{itemize}

Consider now a subdiagram $g^{\half}$ of a Feynman half-diagram $G^{\half}$, with external legs
$\mathbf{z}_{ext}:=\mathbf{z}'_{ext}\uplus \{\xi_v:\ \xi_v\ {\mathrm{uncontracted}}\}.$ The uncontracted
$\phi$-legs $(\xi_v)$ are not considered as {\em true}, free external legs since they are attached on
the mirror and must eventually be integrated, see e.g. eq. (\ref{eq:2.6,5}). Hence one sets \BEQ
\tau_{g^{\half}}A_{g^{\half}}(\mathbf{z}_{ext}):=A_{g^{\half}}(\mathbf{z}'_{ext}=\mathbf{0}, \{\xi_v:\ \xi_v\
{\mathrm{uncontracted}}\}).\EEQ For this reason, it is more natural to write $\tau_g A_{g^{\half}}$
instead of $\tau_{g^{\half}}A_{g^{\half}}$, where in the symmetric graph $g:=(g^{\half})^2$, the
uncontracted $\phi$-legs have now become {\em internal legs}.


\subsection{Diverging graphs}


Consider a connected Feynman diagram $G$. In order to decide whether to renormalize it or not, we
compute its degree of divergence $\omega(G)$. It is simply obtained as the sum of the overall degree of
homogeneity of the integrand, $(1-2\alpha)I_{\phi}(G)-I_{\sigma}(G)$, and of the number,
$I(G)-|V(G)|+1$, of  independent momenta, with respect to which the integrand is integrated; hence it is
simply the overall homogeneity degree of the Feynman integral. Taking into account the relation $|V(G)|=
2I_{\phi}(G)+N_{\phi}(G)$ (obtained by counting one half double line per vertex,
 except for external
double lines which are only connected to one vertex), yields \BEQ \omega(G)=
1-\alpha|V(G)|-(1-\alpha)N_{\phi}(G).\EEQ

\begin{Definition}[diverging graphs] \label{def:div}

We call a  Feynman graph $G$ {\em diverging} if and only if it has no external $\phi$-legs.

\end{Definition}

Clearly enough, with this definition, small graphs (i.e. with  $\alpha |V(G)|<1$) are diverging if and
only if $\omega(G)>0$ (which is the usual definition). It is natural to extend this notion to Feynman
{\em half-diagrams} by letting $\omega(G^{\half}):=1-
\alpha|V(G^{\half})|-(1-\alpha)N_{\phi}(G^{\half})$, where $N_{\phi}(G^{\half})$ is the number of {\em
true} external $\phi$-legs. Then Feynman half-diagrams $G^{\half}(\T)$
 associated to skeleton integrals
of order $n=|V(\T)|<\lfloor 1/\alpha\rfloor$ are diverging if and only if $\omega(G^{\half})>0$.

Consider a connected half-diagram $g^{\half}\subset G^{\half}$ and  its symmetric double
$g:=(g^{\half})^2$. If $g^{\half}$ has uncontracted $\phi$-legs, then $g^{\half}$ is connected to its
image in the mirror by some "bridge", hence $g$ is also connected; $g$ is then called a {\em bilateral diagram}. Otherwise, $g$ is made up of two {\em
unilateral (full) diagrams}. As we shall see in section 4, renormalizing $g^{\half}$ if $g^{\half}$ is
divergent amounts to replacing $\omega(g^{\half})$ with $\omega^*(g^{\half})=\omega(g^{\half})-1$. On
the other hand, if $g^{\half}$ is convergent, no
 renormalization   is performed, hence simply $\omega^*(g^{\half})=\omega(g^{\half})$. Now
 power-counting
 must really be understood in terms of {\em  half-diagrams}, which implies the following rules:

 -- if $g^{\half}$ is a unilateral diagram, then $\omega^*(g^{\half})=\omega(g^{\half})$ if $g^{\half}$
 has external $\phi$-legs, $\omega(g^{\half})-1$ otherwise. Hence (considering that any non-empty
 diagram
 contains at least one line and two vertices) $\omega^*(g^{\half})\le -\alpha$ in any
 case;

 -- if $g^{\half}$ is connected to its image by some bridge, then $\omega^*(g)=\omega(g)$ if $g$
 has external $\phi$-legs, $\omega(g)-2$ otherwise. But since $g$ is symmetric,
 $N_{\phi}(g)$ is {\em even} and $|V(g)|\ge 4$. Hence, in both cases, $\omega(g)\le -1-2\alpha$.

\medskip

These elementary power-counting arguments are essential for section 4.


\subsection{The multiscale BPHZ algorithm}


We denote hereafter by ${\cal F}^{div}(G^{\half})$ the set of forests of diverging subgraphs of
$G^{\half}$. Equivalently, \BEQ {\cal F}^{div}(G):=\{g=(g^{\half})^2\ | \ g^{\half}\in {\cal
F}^{div}(G^{\half})\}, \EEQ if $G=(G^{\half})^2$ is a symmetric graph, is the set of diverging {\em
symmetric} subgraphs.

 We refer to \cite{Riv} or \cite{FVT} for the whole paragraph.
 \medskip

\begin{Definition}[Bogolioubov's non-recursive definition of renormalization]  \label{def:3.2}

\begin{itemize} \item[(i)]

Let

\BEQ {\cal R}A_{G^{\half}(\T)}((\zeta_1,0),\mathbf{\xi}):=\sum_{\F\in{\cal F}^{div}(G(\T))} \prod_{g\in
\F} (-\tau_g)A_{G^{\half}(\T)} ((\zeta_1,0),\mathbf{\xi})). \label{eq:3.2} \EEQ

\item[(ii)] Define correspondingly, for $\nu:=D(f)\mu= {\cal F}^{-1}(f\cdot ({\cal F}\mu))$, with
    $\mu=\otimes_{v\in V(\T)} dB_{x_v}(\ell(v))$  and $f=f(\xi_1,\ldots,\xi_n)$ such that
    $\supp(f)\subset \R^{\T}_+$, \BEA &&  \phi^t_{\nu}(\T):= (2\pi c_{\alpha})^{-n/2}   \int
    \prod_{v\in V(\T)} dW_{\xi_v}(\ell(v))  \frac{e^{\II t\zeta_1}}{[\II\zeta_1]}  d\zeta_1
    \nonumber\\ &&\qquad f(\mathbf{\xi}) {\cal
    R}A_{G^{\half}(\T)}(\mathbf{\zeta}_{ext}=(\zeta_1,0),\mathbf{\xi}_{ext}=(\xi_v)_{v\in V(\T)} ) ,
    \label{eq:3.5} \EEA so that, assuming all decorations $(\ell(v))_{v\in V(\T)}$ are distinct,
    \BEA &&   \Var \left(\phi^t_{\nu}(\T)-\phi^s_{\nu}(\T) \right)=(2\pi c_{\alpha})^{-n} \int
    \frac{d\zeta_1}{\zeta_1^2} |e^{\II t\zeta_1}-e^{\II s\zeta_1}|^2 \Var
    \hat{\phi}^{\zeta_1}_{\nu}(\T),\nonumber\\ && \qquad   \Var
    \hat{\phi}^{\zeta_1}_{\nu}(\T)=D(f){\cal R} A_{G(\T)}(\mathbf{\zeta}_{ext} =(\zeta_1,0)), \EEA
    where \BEQ D(f) {\cal R}A_{G(\T)}(\zeta_1,0):=  \int \prod_{v\in V(\T)} d\xi_v  f^2(\mathbf{\xi})
 \left|{\cal R}A_{G^{\half}}((\zeta_1,0),\mathbf{\xi})\right|^2.
\EEQ

\end{itemize} \end{Definition}

Now come two  essential remarks, based on the fact  that divergent subgraphs have {\em no} external
$\phi$-leg by definition.

\begin{enumerate} \item Since renormalization leaves $\xi$-momenta unchanged, one may consider the
integration measure $f(\mathbf{\xi})\prod_{v\in V(\T)}
    dW_{\xi_v}(\ell(v))$ in eq. (\ref{eq:3.5}) as a simple decoration of the
    vertices. In this sense $\phi^t_{\nu}(\T)$ may be considered as a {\em renormalized skeleton
    integral}, denoted by
 ${\cal R}\overline{\SkI}^t_{\nu}(\T)$.

\item
    Consider some multiple contraction
$\phi^t_{\nu}(\T;(i_1 i_2),\ldots,(i_{2p-1},i_{2p}))$ of $\phi^t_{\nu}(\T)$. Then \BEA &&
\phi^t_{\nu}(\T;(i_1 i_2),\ldots,(i_{2p-1} i_{2p})):= (2\pi c_{\alpha})^{-n/2}    \int \prod_{v\in
V(\T)} dW_{\xi_v}(\ell(v)) \frac{e^{\II t\zeta_1}}{[\II\zeta_1]} d\zeta_1 \nonumber\\
 &&\qquad f(\mathbf{\xi})  {\cal R}A_{G^{\half}(\T;(i_1 i_2),\ldots,(i_{2p-1}
 i_{2p}))}(\mathbf{\zeta}_{ext}=(\zeta_1,0),\mathbf{\xi}_{ext}=(\xi_v)_{v\in
V(\T)} ).  \nonumber\\ \EEA In other words, contractions and renormalization {\em commute}. This
remark extends in a straightforward way to contractions between different trees as in Lemma
\ref{lem:int-Feynman} (2). This allows us to extend the BPHZ construction to contracted graphs.
Namely, consider the Feynman diagram $G=(G^{\half})^2$ obtained by gluing two identical Feynman
half-diagrams {\em with the same external structure}, i.e. such that $\bar{z}=z$ whenever $z$ is a
{\em true} external leg. Then {\em all} (internal or external) momenta $\zeta$ or $\xi$ are equal to
their image $\bar{\zeta}$ or $\bar{\xi}$ in the mirror. Now one defines \BEQ {\cal
R}A_G(z_{ext})=\int \prod_{\xi\ |\ \xi \ {\mathrm{uncontracted}}} d\xi |{\cal
R}A_{G^{\half}}(z_{ext},\xi)|^2 \label{eq:twice-ren} \EEQ where ${\cal R}A_{G^{\half}}( \ .\
)=\sum_{\F\in {\cal F}^{div}(G)} \prod_{g\in\F} (-\tau_g) A_{G^{\half}} (\ .\ )$ is defined by the
BPHZ formula as in eq. (\ref{eq:3.2}).

\item Let $G=(G^{\half})^2$ be as in 2. As already mentioned, ${\cal R}A_G(\ .\ )$ differs from the
    usual BPHZ renormalized graph amplitude since (due to the square in the right-hand side of eq.
    (\ref{eq:twice-ren})) divergent bilateral subgraphs are in some sense renormalized {\em twice} from the
    point of view of power-counting.

\end{enumerate}

\bigskip

Choose some constant $M>1$. An {\em attribution of momenta} $\mu$ for a Feynman diagram $G$ is a choice
of {\em $M$-adic scale} for each momentum of $G$, i.e. a function $\mu:L(G)\cup L_{ext}(G)\to \Z$ and an
associated restriction of the momentum $|z_{\ell}|=|\zeta_{\ell}|$ or  $|\xi_{\ell}|$, $\ell\in L(G)\cup
L_{ext}(G)$ to the $M$-adic interval $[M^{\mu(\ell)},M^{\mu(\ell)+1})$. Thus one may define, e.g. for  a tree
Feynman half-diagram $G^{\half}$, compare with eq. (\ref{eq:2.1}), \BEA &&
A_{G^{\half}}^{\mu}(\mathbf{\zeta}_{ext},\mathbf{\xi}_{ext}):= \del(\zeta_{ext}+\xi_{ext})\nonumber\\ && \quad
\prod_{v\in V(\T)\ |\ \xi_v \ {\mathrm{uncontracted}}}
 |\xi_v|^{\half-\alpha}. \int \prod_{q=1}^p |\xi_{(i_{2q-1} i_{2q})}|^{1-2\alpha} d\xi_{(i_{2q-1}
 i_{2q})} \nonumber\\
 &&\qquad
  \prod_{v\in V(\T)\setminus\{ {\mathrm{roots}} \}} \frac{1}{\zeta_v}
 \prod_{v\in V(\T)} \left(  {\bf 1}_{|\zeta_v|\in  [M^{\mu(\zeta_v)},M^{\mu(\zeta_v)+1} )} \right)
 \left(  {\bf 1}_{|\xi_v|\in  [M^{\mu(\xi_v)},M^{\mu(\xi_v)+1} )} \right)
 \nonumber\\ \EEA
where by definition $M^{\mu(\xi_{i_{2q-1}})}=M^{\mu(\xi_{i_{2q}})}=M^{\mu(\xi_{(i_{2q-1} i_{2q})})}$ for
contracted lines, and similarly for an arbitrary Feynman diagram $G$, compare with eq. (\ref{eq:2.3}),
\BEA &&  A^{\mu}_G(\mathbf{z}_{ext}):=\del(z_{ext}) \int \prod_{\ell\in L(G)\setminus L'(G)} dz_{\ell}
\prod_{\ell\in L(G)\cup L_{ext}(G)} {\bf 1}_{z_{\ell}\in [M^{\mu(\ell)},M^{\mu(\ell)+1})}  \nonumber\\
&& \qquad \prod_{\ell\in L_{\phi}(G)} |\xi_{\ell}|^{1-2\alpha} \prod_{\ell\in L_{\sigma}(G)}
\zeta_{\ell}^{-1}.\EEA

Feynman diagrams with a  fixed scale attribution are called {\em multiscale diagrams}. In the corresponding graphical representation (see below), vertices are split according to the scales of the lines attached to them.

\begin{Definition}[Gallavotti-Nicol\`o tree]  \label{def:GN}

Let $G^{j}\subset G$, $j\in\Z$ be the subdiagram with set of lines $L(G^{j})\cup L_{ext}(G^j):=
\{\ell\in L(G)\cup L_{ext}(G);\ \mu(\ell)\ge j\}$, and $(G^{j}_k)_{k=1,2,\ldots}$ the connected
components of $G^{j}$.

Thet set of connected subgraphs $(G^{j}_k)_{j,k}$ -- called {\em local subgraphs} --  makes up a tree of subgraphs of $G$, called {\em
Gallavotti-Nicol\`o tree}.

\end{Definition}

\begin{figure}[h]
  \centering
   \includegraphics[scale=0.35]{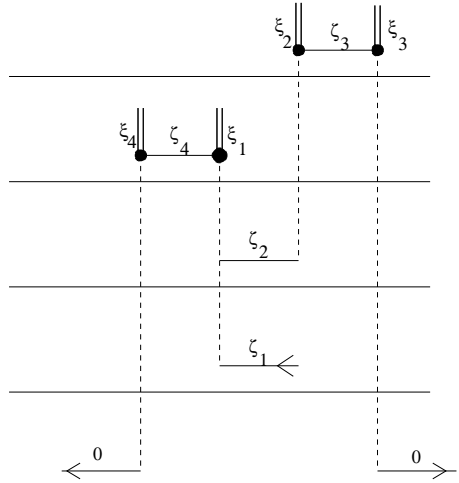}
   \caption{\small{A Gallavotti-Nicol\`o tree (case 1).}}
  \label{Fig-renbis3}
\end{figure}

\begin{figure}[h]
  \centering
   \includegraphics[scale=0.35]{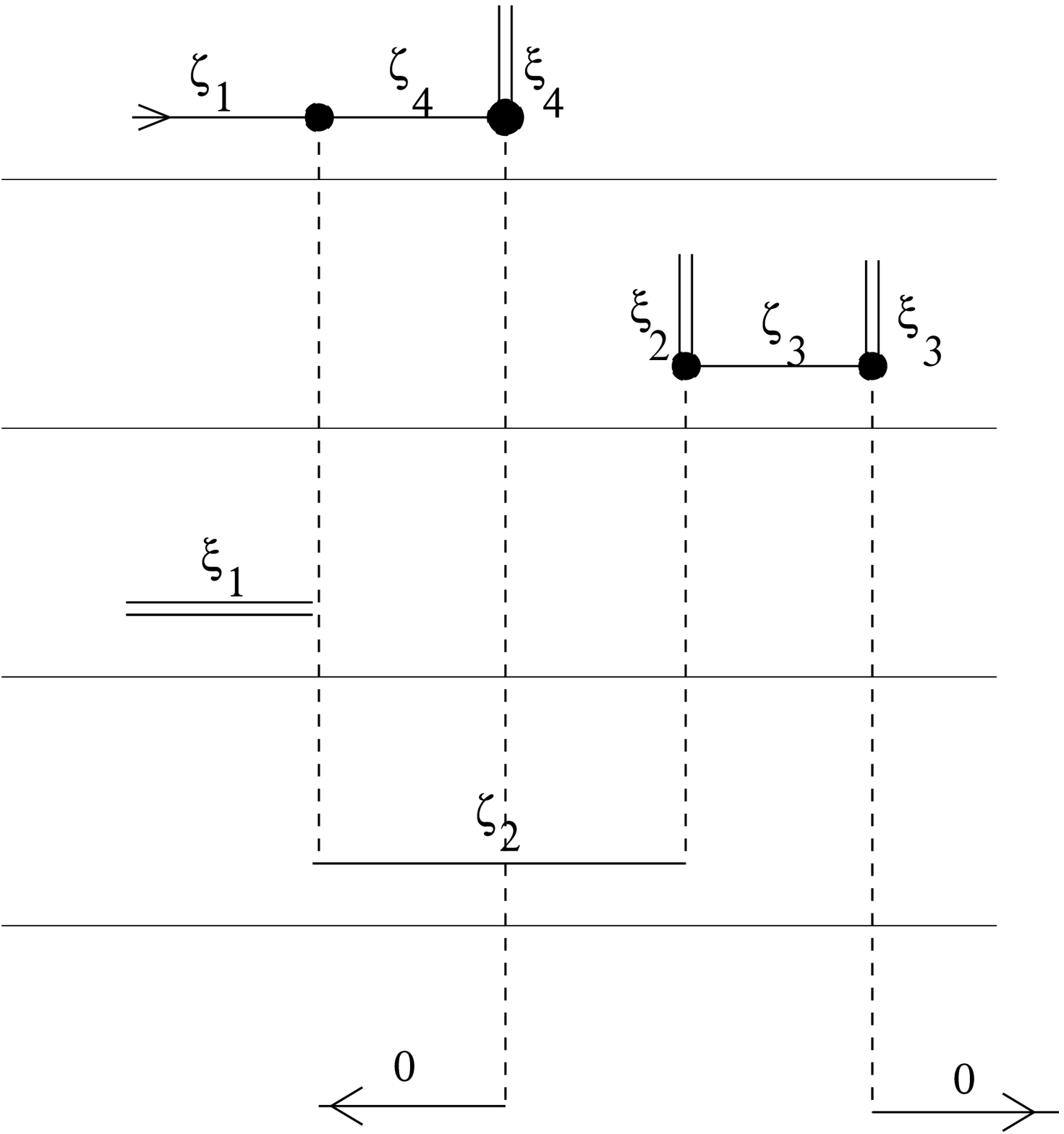}
   \caption{\small{A Gallavotti-Nicol\`o tree (case 2).}}
  \label{Fig-renbis4}
\end{figure}

Two instances of Gallavotti-Nicol\`o trees are represented on Fig. \ref{Fig-renbis3}, \ref{Fig-renbis4}.
By shifting slightly the $M$-adic intervals, it is possible to manage to have both lines of highest
momentum of any given vertex in the same interval.

\begin{Definition}

 Let $\F\in {\cal F}^{div}(G)$ be a forest of diverging  subgraphs of $G$.

\begin{itemize} \item[(i)] Let $g\in G$ be a subgraph of $G$. Then $g$ is {\em compatible with} $\F$ if
and only if $\F\cup \{g\}$ is a forest. \item[(ii)] Assume $g\in G$ is compatible with $\F$. We let
$g^-_{\F}$ be the ancestor of $g$ in the forest of graphs $\F\cup\{G\}$, and $g^{\uparrow}_{\F}$ be the
union of its children, namely, \BEQ g^{\uparrow}_{\F}=\cup_{h\subsetneq g,h\in \F} h.\EEQ \item[(iii)]
Let $\mu$ be a momentum scale attribution.  The {\em dangerous forest} $D^{\mu}(\F)\subset\F$
 associated
to the forest $\F$ and the momentum scale attribution $\mu$ is the sub-forest defined by \BEQ
\left(g\in D^{\mu}(\F)\right)\Longleftrightarrow \left( \min\{i_{\ell}(\mu)\ :\ \ell\in L(g\setminus
g^{\uparrow}_{\F})\} > \max\{i_{\ell}(\mu) \ :\ \ell\in L_{ext}(g)\cap L(g^-_{\F})\} \right).\EEQ
\item[(iv)]  Call the sub-forest $ND^{\mu}(\F):=\F\setminus D^{\mu}(\F)\subset\F$ the {\em
non-dangerous} or {\em harmless forest} associated to $\F$ and $\mu$. \end{itemize}

\end{Definition}

One can prove that $ND^{\mu}\circ ND^{\mu}=ND^{\mu}$. Hence \BEQ {\cal F}^{div}(G)=\cup_{\F\in {\cal
F}^{div}(G) \ |\ T_{\mu}(\F)=\F} \{\F'\supset \F\ :\ ND_{\mu}(\F')=\F\}. \EEQ One obtains the following
classification of forests:

\begin{Proposition}

Let

\begin{itemize} \item[(i)]
 $Safe^{\mu}(G)\subset {\cal F}^{div}(G)$ be the set of forests of diverging graphs which are
 invariant under
the projection operator $ND^{\mu}$ and thus harmless,
 namely, $Safe^{\mu}(G):=\{\F\in {\cal F}^{div}(G)\ :\ ND^{\mu}(\F)=\F\}$;
\item[(ii)] $Ext^{\mu}(\F)\subset {\cal F}^{div}(G)$, with $\F\in Safe^{\mu}(G)$, be the ``maximal
    dangerous extension'' of the harmless forest $\F$ within the $ND^{\mu}$-equivalence class of
    $\F$, namely, $\F\uplus Ext^{\mu}(\F)$ is the maximal forest such that $ND^{\mu}(\F\uplus
    Ext^{\mu}(\F))=\F$. \end{itemize}

Then: \begin{itemize} \item[(i)] \BEQ \left( ND^{\mu}(\F')=\F \right)\Longleftrightarrow \left(
\F\subset \F'\subset \F\uplus Ext^{\mu}(\F) \right);\EEQ \item[(ii)] $Ext^{\mu}(\F)$ is the set of
subgraphs $g\in G$, compatible with $\F$, such that $g\in D^{\mu}(\F\cup\{g\})$. \end{itemize}

\end{Proposition}

In particular, $Ext^{\mu}(\emptyset)$ is the forest of local subgraphs of $G$, or in other words
the Gallavotti-Nicol\`o tree, see Definition \ref{def:GN}.

\begin{Corollary}

\BEQ {\cal R}A_{G^{\half}}=\sum_{\F\in {\cal F}^{div}(G)} {\cal R}A_{G^{\half},\F} \EEQ where \BEQ {\cal
R}A_{G^{\half},\F}:=\sum_{\mu\ |\ \F\in Safe^{\mu}(G)} \prod_{g\in \F} (-\tau_g) \prod_{h\in
Ext^{\mu}(\F)} (1-\tau_h)A^{\mu}_{G^{\half}}.\EEQ

\end{Corollary}

The BPHZ renormalization scheme is perfect in perturbative field theory, but experts of constructive field theory scorn it because it leads to unwanted combinatorial factors of order 
$O(n!)$ for large Feynman diagrams with $O(n)$ vertices, called renormalons (see \cite{FVT}, eq.
(1.1.12)), which ruin any hope of resumming the series of perturbations. These may be avoided
by considering only {\em useful renormalizations} associated to {\em local subgraphs} in the 
sense of Definition \ref{def:GN}, at the price of introducing scale-dependent renormalized
coupling constants, see \cite{FVT}, \S 1.4. This gives another possible renormalization formula,
which is however scale-dependent,
\BEQ {\cal R}^{useful}A_{G^{\half}}={\cal R}A_{G^{\half},\emptyset}=\sum_{\mu}\prod_{j,k}
(1-\tau_{G^j_k})A^{\mu}_{G^{\half}}.\EEQ

In our context, the whole discussion seems a priori pointless since (i) required Feynman diagrams have
at most $2\lfloor 1/\alpha\rfloor<\infty$ vertices; (ii) there are no coupling constants at all.
Purely esthetic reasons plead for the scale-independent renormalization ${\cal R}A_{G^{\half}}$.
However, it may be that using ${\cal R}^{useful}{A_{G^{\half}}}$ instead of ${\cal R}A_{G^{\half}}$ gives better bounds for higher-order iterated integrals, which may after all
also be rewritten as Feyman diagrams. Good bounds are notoriously difficult to obtain for
general rough paths, which is a major problem when solving stochastic differential equations,
see \cite{FV} for a general discussion, or \cite{Unt-lin} in the particular case of linear
stochastic differential equations, in connection with the Magnus series.


\section{Main bound for Feynman diagrams}


This section is devoted to the proof by classical multi-scale arguments \cite{Riv,FVT}
 of the following theorem.

\begin{Definition}[highest bridge] (see end of subsection 3.1)

Let $g=(g^{\half})^2$ be a multi-scale symmetric Feynman diagram, such that $g^{\half}$ is connected. If
$g^{\half}$ has at least one uncontracted $\xi$-leg, then $g$ is connected by "bridges". Then the {\em
highest bridge} is the uncontracted $\xi$-leg of highest scale.

\end{Definition}

Let $G=G(\T;(i_1 i_2)\ldots(i_{2p-1} i_{2p}))$
 be a  symmetric tree Feynman diagram with $2n$ vertices: then
 $G=(G^{\half})^2$ is made up of two disconnected unilateral Feyman diagrams if and only if
 $G$ has been totally contracted, i.e. $2p=n$, in which case the momentum conservation condition implies
 that $\zeta_{ext}=0$.  Then there is no bridge and hence no highest bridge.  In particular, if $\T$ is connected, so that $\mathbf{\zeta}_{ext}=\{\zeta_1\}$, the diagram evaluation $A_G$ vanishes by symmetry (namely, $\zeta_1=0$, and the denominator $\frac{1}{\zeta_2\ldots\zeta_{2n}}$ changes sign when all momenta are changed to their opposites). On the other hand, assuming
 $n<\lfloor 1/\alpha\rfloor$,  $\omega(G^{\half})=1-n\alpha>0$, whereas $\omega(G)=1-2n\alpha<0$ if
 $\frac{1}{2\alpha}<n<\frac{1}{\alpha}$ for a {\em connected}, symmetric tree Feynman diagram.

\medskip

Estimates for Feynman diagrams with $2n$ vertices must be expressed in terms of a {\em reference scale}. It turns out that any (internal or external) momentum may be chosen as a reference scale when $2n<\lfloor
1/\alpha \rfloor$, because the renormalized amplitude is then both ultra-violet and infra-red convergent.
On the other hand, diagrams with $1/\alpha<2n<2/\alpha$ vertices (thus not unilateral) increase indefinitely when external momenta go to zero, and computations show that momenta above the highest bridge are too "loosely" attached to those below to control the infra-red behaviour of the whole diagram. In that case, the most appropriate reference scale is that of the highest bridge. This is the
content of the following Theorem.

\bigskip

\begin{Theorem} \label{th:4.1}

Let $G:=G(\T;(i_1 i_2)\ldots (i_{2p-1}i_{2p}))$ be a symmetric tree Feynman
 diagram with $2n<2/\alpha$ vertices. Write  $\mathbf{\zeta}_{ext}=(\zeta_{r_1},\ldots,
 \zeta_{r_q},\bar{\zeta}_{r_1},\ldots,\bar{\zeta}_{r_q})$ as in Lemma \ref{lem:int-Feynman}. Assume $\zeta_{r_m}=\bar{\zeta}_{r_m}$,
  $m=1,\ldots,q$, so that each $\zeta$-momentum and each contracted  $\xi$-momentum is
  equal to the corresponding $\bar{\zeta}$- or $\bar{\xi}$-momentum on the other side of the mirror.

  Label $\mathbf{\zeta}_{ext}$ so that $|\zeta_{r_1}|<\ldots<|\zeta_{r_q}|$.

\begin{enumerate}

\item (bilateral diagrams)

 Assume $G$ is bilateral, so $G$ is connected. Let
$\xi_{ref}$ be the highest bridge.   {\em Fix} $j_{ref}:=j(\xi_{ref})$ and sum over all scale
 attributions $\mu$ such that $\mu(\xi_{ref})=j_{ref}$. Replace one of the $\xi$-propagators, $|\xi_1|^{1-2\alpha}$, say, by $|\xi_1|^{(1-2\alpha)-2n'\alpha}$ in the integrand, with
 $n'\ge 0$, $n+n'<1/\alpha$. Denote by
 ${\cal R}A_{G \stackrel{\rightarrow}{\xi_1} n'}^{j_{ref}}(\mathbf{\zeta}_{ext}):=\sum_{\mu} {\cal R}A_{G  \stackrel{\rightarrow}{\xi_1} n'}^{\mu}(\mathbf{\zeta}_{ext})$ the
  result. Then:

 \BEQ \Var {\cal R}A_{G \stackrel{\rightarrow}{\xi_1} n'}^{j_{ref}}(\mathbf{\zeta}_{ext})\lesssim M^{(1-2(n+n')\alpha)j_{ref}}
   \left(\frac{\min(|\zeta_{r_1}|,M^{j_{ref}})}{\max(|\zeta_{r_q}|,M^{j_{ref}})} \right)^{\alpha^-}
 \label{eq:4.4} \EEQ
 whenever $\alpha^-<\alpha$.

\item (diagrams with $2n<\lfloor 1/\alpha\rfloor$ vertices)

Let $G$ be indifferently a unilateral diagram, or a bilateral  with
$2n<\lfloor 1/\alpha\rfloor$ vertices. Let $\xi_{ref}$ be one of the $\xi$-lines of $G$. Fix $j_{ref}:=j(\xi_{ref})$ and sum  over all scale attributions $\mu$ such that $\mu(\xi_{ref})=j_{ref}$. Replace the propagator  $|\xi_{ref}|^{1-2\alpha}$, say, by $|\xi_{ref}|^{(1-2\alpha)-2n'\alpha}$ in the integrand, with
 $n'\ge 0$, $n+n'<1/\alpha$.  Denote by
 ${\cal R}A_{G \stackrel{\rightarrow}{\xi_{ref}} n'}^{j_{ref}}(\mathbf{\zeta}_{ext})$ the result. Then eq. (\ref{eq:4.4}) holds.

 \end{enumerate}

\end{Theorem}

{\bf Remarks.}

 \begin{enumerate}

 \item The factor $\left(\frac{\min(|\zeta_{r_1}|,M^{j_{ref}})}{\max(|\zeta_{r_q}|,M^{j_{ref}})}
 \right)^{\alpha^-}$
 in eq. (\ref{eq:4.4}) is obtained
and shall be used as a product of {\em spring factors}, $\prod_{m=1}^{q}
\left|\frac{u_m}{u_{m+1}}\right|^{\alpha^-}$, where $|u_1|<\ldots<|u_{q+1}|$ is the ordered list of
momenta $(M^{j_{ref}},|\zeta_{r_1}|,\ldots,|\zeta_{r_q}|)$.

\item The supplementary factors $|\xi_1|^{-2n'\alpha}$ or $|\xi_{ref}|^{-2n'\alpha}$ may be seen as a "grafting" of another tree $\T'$ on $\T$. It will be used for $G=G_1$ and unrooted diagrams $G'_i,i=1,\ldots,I'$ (see introduction to section 5). The term "grafting" is only approximate since $\T$ and $\T'$ remain disjoint.

\end{enumerate}

\medskip

Using the Cauchy-Schwarz inequality in eq. (\ref{eq:2.2}), this result yields immediately

\begin{Corollary}

Consider a bilateral diagram $G$. Then

 \BEQ \Var {\cal R}A_{G\stackrel{\to}{\xi_1}n'}^{j_{ref}}(\mathbf{\zeta}_{ext},\mathbf{\bar{\zeta}}_{ext})\lesssim
 M^{(1-2(n+n')\alpha)j_{ref}}
   \left(\frac{\min(|\zeta_{r_1}|,M^{j_{ref}})}{\max(|\zeta_{r_q}|,M^{j_{ref}})} \right)^{\alpha^-/2}
   \left(\frac{\min(|\bar{\zeta}_{r_1}|,M^{j_{ref}})}{\max(|\tilde{\zeta}_{r_q}|,M^{j_{ref}})}
   \right)^{\alpha^-/2}.
 \label{eq:4.4bis} \EEQ

\end{Corollary}

{\bf Proof of Theorem \ref{th:4.1}.}

Let $\mu$ be any attribution of momenta. We shall consider only {\em useful} renormalizations
in the proof. Namely, as shown in \cite{FVT}, \S 1.3, the operations $\prod_{g\in\F} (-\tau_g)$,
$\F\in Safe^{\mu}(G)$, see eq. (), are equivalent to displacing all external $\zeta$-legs to
the same point,  and do not change the power-counting rules in the proof.

Choose inductively, starting from the highest momentum scale, a subset of lines $L'(G)\subset L(G)$ so
that $(z_{\ell})_{\ell \in L(G^{j}_k)\setminus L'(G^{j}_k)}$, where $L'(G^{j}_k):=L'(G)\cap L(G^{j}_k)$,
make up a maximal set of independent momenta of the graph $G^{j}_k$ shorn of its external legs. Note that $G^j_k$ is necessarily symmetric. The degree of divergence $\omega(G^j_k)=1-2|V(G^j_k)|\alpha-(1-\alpha)N_{\phi}(G^j_k)$ has been defined in subsection 3.1.

\medskip

Assume for a moment that all $(z_{\ell})_{\ell\in L(G^{j}_k)}$ are of the same order, $M^{j'}$, say
$(j'\ge j)$. Then the previous power-counting arguments show that $A_{G^j_k}$ is of order
$M^{j'\omega(G^j_k)}$. If $\omega(G^{j}_k)\ge 0$, then clearly the sum $\sum_{j'=j}^{\infty}
M^{j'\omega(G^{j}_k)}$ diverges, so the sum over all momenta attributions diverges. On the other hand, if
$\omega(G^{j}_k)<0$, then the sum over all momenta attributions may still diverge because of so-called
{\em sub-divergences} due to  the higher subgraphs $G^{j'}_k$, $j'>j$; the graph is a priori only
 {\em overall
convergent}.

\medskip

Let us now see how renormalization will make all symmetric subgraphs convergent. Consider any of the local
subgraphs $G^j_k$. Assume $\omega(G^j_k)>0$, so that $G^j_k$ must be renormalized. We introduce some
notations for the sake of clarity.  Let quite generally $V_{ext}(g)$ be the set of external vertices of a graph $g$, and $L_v$, resp. $L_{v,ext}$ ($v\in V(G)$) be the set of internal, resp. external lines
of $g$ attached to $v$. Now, to each $v\in V_{ext}(G^j_k)$, one associates the unique line
$\ell'_v \in L'(G^j_k)\cap L_v(G^j_k)$, and lets $z_v^*:=\sum_{\ell\in L_v(G^j_k)\setminus \{\ell'_v\}}
z_{\ell}$ and $z_{v,ext}:=\sum_{\ell \in L_{v,ext}(G^j_k)} z_{\ell}$, so that
$z_{\ell'_v}+z_v^*+z_{v,ext}=0$. Choose some arbitrary ordering of the external legs of $G$,
$L_{ext}(G^j_k)=\{\ell_1,\ldots,\ell_{|L_{ext}(G^j_k)|}\}$.  Renormalization changes only the values
of the external momenta, so it acts really on the product $A_{ext}(G^j_k):=\prod_{v\in V_{ext}(G^j_k)}
|z_v^*+z_{v,ext}|^{\beta_{\ell'_v}}$, with $\beta_{\ell'_v}=1-2\alpha$ or $-1$,

 \BEA  &&  A_{ext}(G^j_k)\rightsquigarrow  {\cal
R}A_{ext}(G^j_k):=\prod_{v\in V_{ext}(G^j_k)} |z_v^*+z_{v,ext}|^{\beta_{\ell'_v}} - \prod_{v\in
V_{ext}(G^j_k)} |z_v^*|^{\beta_{\ell'_v}} \nonumber\\ \qquad \qquad &&=\sum_{i=1}^{|L_{ext}(G^j_k)|}
\int_0^1 z_{\ell_i} \partial_{z_{\ell_i}} \prod_{v\in V_{ext}(G^j_k)}
|z_v^*+sz_{v,ext}|^{\beta_{\ell'_v}} ds. \EEA

Now only one or two factors in the product over $V_{ext}(G^j_k)$ depend on $z_{\ell_i}$; the derivative
$\partial_{z_{\ell_i}}$ acting on each of these, generically denoted by
$|z_v^*+sz_{v,ext}|^{\beta_{\ell'_v}}$, generates an extra multiplicative factor called  {\em
spring factor},
  $\frac{z_{\ell_i}}{z_v^*+sz_{v,ext}}$, up to a constant.  This spring factor is at most $O(M^{\min\mu(G^j_k)-\max\mu(\partial G^j_k)})$, where
$\min\mu(G^j_k)$ is the minimal scale index of all {\em internal} lines of $G^j_k$, and
$\max\mu(\partial G^j_k)$ the maximal scale index of all {\em external} lines of $G^j_k$.

We shall now rewrite ${\cal  R}A(G)$ by using the local graph decomposition of $G$. First, each factor
$M^{\beta_{\ell}}$, $\ell\in L(G)$ may be rewritten as $\prod_{(j,k) ; \ell\in L(G^j_k)}
M^{\beta_{\ell}}$. Similarly, the integration over the independent momenta yields $\prod_{(j,k)}
M^{|L(G^j_k)|-|V(G^j_k)|+1}$.  Multiplying these two expressions, one gets $\prod_{(j;k)}
M^{\omega(G^j_k)}$. Now,  the spring factors due to renormalization contribute -- due
 to the fact that $A_G$ and ${\cal R}A_G$ are {\em squared amplitudes} -- a factor
$M^{-2}$ per scale until $G^j_k$ absorbs one external line, so all together
$\prod_{(j,k);\omega(G^j_k)>0} M^{-2}$, where unilateral diagrams are only counted {\em once}.
 Finally, "grafting" $|\xi_1|^{-2n'\alpha}$ into the graph
  is equivalent to subtracting $2n'\alpha$ to all $\omega(G^j)$ with $j\le j(\xi_1)$, with
  $\xi_1=\xi_{ref}$ in case (2). All together, one has proved that
 \BEQ {\cal R}A^{\mu}(G\stackrel{\rightarrow}{\xi_1} n')\le K^{2n} \prod_{(j,k)} M^{\omega^*(G^j_k)} \ \cdot\ \prod_{j\le j(\xi_1)} M^{-2n'\alpha} \EEQ
for some constant $K$, where $\omega^*(G^j_k)=\omega(G^j_k)$ if
$\omega(G^j_k)<0$, and $\omega(G^j_k)-1$, resp. $\omega(G^j_k)-2$ otherwise for unilateral, resp.
bilateral subdiagrams, {\em except} for the total graph $G$ which is not renormalized (having only
external legs of zero momentum), so that $\omega^*(G)=\omega(G)=1-2n\alpha$. As noted at the end of
subsection 3.1, $\omega^*(G^j_k)\le -\alpha$, resp. $\le -1-\alpha$ for unilateral, resp. bilateral
subdiagrams others than $G$. Summing up the divergence degrees of a given scale $j$,
$\omega^*(G^{j})=\sum_k \omega^*(G^{j}_k)$, yields a quantity $\le -1-\alpha$ if $j\le j_{ref}$
in case 1  because
the subdiagram containing the highest bridge is bilateral. Finally,
${\cal R}A^{\mu}(G\stackrel{\rightarrow}{\xi_1}n')\le K^{2n} \prod_j M^{\omega^*_{gr}(G^j)}$ if
one lets $\omega^*_{gr}(G^j):=\omega^*(G^j)-2n'\alpha$ ($j\le j(\xi_1)$), $\omega^*(G^j)$
($j>j(\xi_1)$) be the equivalent degree of divergence of $G^j$ after renormalization and grafting.

\medskip

Fix the scales of $\mu$, say, $j_1<j_2<\ldots,j_I=j_{max}$, with $j_{I_1}=j$, and let $j'_1<\ldots<
j'_{q+1}$ be the scales of $|\zeta_{r_1}|,\ldots,|\zeta_{r_q}|,M^{j_{ref}}$ put into increasing order.
 Then the renormalized amplitude is bounded up to a constant  by
\BEA&&  \sum_{j_1>-\infty} M^{j_1 \omega^*_{gr}(G)} \left( \sum_{j_2\ge j_1} M^{(j_2-j_1)\omega^*_{gr}(G^{j_2})}
 \left(
\cdots \left( \sum_{j_I \ge j_{I-1}} M^{(j_I-j_{I-1})\omega^*_{gr}(G^{j_I})} \right) \cdots \right)\right)
\label{eq:4.5}  \nonumber\\ &&  \le \prod_{m=1}^q M^{-(j'_{m+1}-j'_m)\alpha^-} \ \cdot \
\sum_{j_1>-\infty} M^{j_1 \omega^*_{gr}(G)} \left( \sum_{j_2\ge j_1}
M^{(j_2-j_1)(\omega^*_{gr}(G^{j_2})+\alpha^-)} \right.\nonumber\\ && \qquad \left.
 \left(
\cdots \left( \sum_{j_I \ge j_{I-1}} M^{(j_I-j_{I-1})(\omega^*_{gr}(G^{j_I})+\alpha^-)} \right) \cdots
\right)\right)    \nonumber\\ \EEA the scale $j_{I_1}$ being fixed, and the scales
$j_1,\ldots,j_{I_1-1}$ constrained to be below $j_{ref}$. Since all $\omega^*_{gr}(G^{j_i})$ {\em except}
possibly  $\omega^*_{gr}(G)=\omega(G)-2n'\alpha$ are $\le -\alpha<-\alpha^-$,
 one may sum
down to scale $j_{I_1}$, which  (discarding the $\alpha^-$-spring prefactors) leads to the following bound,

\BEA &&  \sum_{j_1>-\infty} M^{j_1 \omega^*_{gr}(G)} \left( \sum_{j_2\ge j_1}
 M^{(j_2-j_1)(\omega^*_{gr}(G^{j_2})+\alpha^-) }  \right.\nonumber\\
 && \qquad \left.
\left( \ldots \left( \sum_{j_{I_1}\ge j_{I_1-1}} M^{(j_{I_1}-j_{I_1-1})(\omega^*_{gr}(G^{j_{I_1}})+\alpha^-)
}
 \right) \ldots \right)\right). \nonumber\\ \EEA

However, $j_{ref}$ is fixed, hence this expression must be computed as \BEA &&  M^{j_{I_1}
(\omega^*_{gr}(G^{j_{I_1}})+\alpha^-)} \cdot\ \sum_{j_1\le j_{I_1}}
M^{j_1(\omega^*_{gr}(G)-\omega^*_{gr}(G^{j_2})-\alpha^-)} \sum_{j_2=j_1}^{j_{I_1}} M^{j_2(\omega^*_{gr}(G^{j_2})-
\omega^*_{gr}(G^{j_3}))} \ldots \nonumber\\ && \qquad \sum_{j_{I_1-1}=j_{I_1-2}}^{j_{I_1}}
M^{j_{I_1-1}(\omega^*_{gr}(G^{j_{I_1-1}})-\omega^*_{gr}(G^{j_{I_1}}))}, \nonumber\\ \EEA

or (integrating from the lowest to the highest scale instead)

\BEA &&  M^{j_{I_1} (\omega^*_{gr}(G^{j_{I_1}})+\alpha^-)} \cdot\
 \sum_{j_{I_1-1}<j_{I_1}}
M^{j_{I_1-1}(\omega^*_{gr}(G^{j_{I_1-1}})-\omega^*_{gr}(G^{j_{I_1}}))} \nonumber\\ && \ldots
 \sum_{j_2<j_3} M^{j_2(\omega^*_{gr}(G^{j_2})-
\omega^*_{gr}(G^{j_3}))}    \sum_{j_1<j_2} M^{j_1(\omega^*_{gr}(G)-\omega^*_{gr}(G^{j_2})-\alpha^-)}. \label{eq:4.8}
\EEA

This is convergent if and only if  $\omega^*_{gr}(G)-\omega^*_{gr}(G^{j_2})-\alpha^-$,
$\left(\omega^*_{gr}(G)-\omega^*_{gr}(G^{j_2})-\alpha^-\right)+
\left(\omega^*_{gr}(G^{j_2})-\omega^*_{gr}(G^{j_3})
 \right)
=\omega^*_{gr}(G)-\omega^*_{gr}(G^{j_3})-\alpha^-,\ldots$, $\omega^*_{gr}(G)-\omega^*_{gr}(G^{j_{ref}})$ are $>0$. This holds true
since $\omega^*_{gr}(G)-\omega^*_{gr}(G^{j_i})-\alpha^-\ge (1-2(n+n')\alpha)+(1+\alpha)-\alpha^->0$
 in case 1, and $\ge (1-2n\alpha)+\alpha-\alpha^->0$ in case 2. Hence
one gets in the end a bound of order $O(M^{(1-2n\alpha)j_{ref}}).$

\hfill \eop

{\bf Examples.} In the two examples below, we use as reference scale that of the external
$\zeta$-leg, called $\zeta_1$ here by reference to the root of the corresponding tree and
let $n'=0$ to simplify. Taking for reference scale some internal $\xi$-line as in Theorem 
\ref{th:4.1} would of course be possible, with minor differences. 

\begin{enumerate} \item  Consider the first Gallavotti-Nicol\`o tree of Fig. \ref{Fig-renbis3}. One may
choose as integration variables $L(G)\setminus L'(G)=\{\zeta_2,\zeta_3,\zeta_4\}$, so that
$\xi_2=\zeta_2-\zeta_3$, $\xi_3=\zeta_3$, $\xi_4=\zeta_4$, $\xi_1=\zeta_1-\zeta_2-\zeta_4$. Hence \BEQ
A(G)=\int d\zeta_2 d\zeta_3 d\zeta_4  \left(  |\zeta_2-\zeta_3|^{\half-\alpha} |\zeta_3|^{-\half-\alpha}
\ \cdot\ |\zeta_1-\zeta_2-\zeta_4|^{\half-\alpha} |\zeta_4|^{-\half-\alpha}.\ \cdot\ \zeta_2^{-1}
\right)^2. \EEQ

The subdiagrams with lines $(\xi_2,\zeta_3,\xi_3)$, $(\xi_4,\zeta_4,\xi_1)$ are
 renormalized by subtracting their value at $\zeta_2=0$, and then the larger subdiagram
 $(\xi_4,\zeta_4,\xi_1,
\zeta_2,\xi_2,\zeta_3,\xi_3)$ is further renormalized by subtracting its value at $\zeta_1=0$.
Hence
  $|\zeta_2-\zeta_3|^{\half-\alpha}$ is replaced with
  $|\zeta_2-\zeta_3|^{\half-\alpha}-|\zeta_3|^{\half
-\alpha}=O(\zeta_2\ \cdot |\zeta_3|^{-\half-\alpha})$, and
$|\zeta_1-\zeta_2-\zeta_4|^{\half-\alpha}$ by \BEQ
\left(|\zeta_1-\zeta_2-\zeta_4|^{\half-\alpha}-|\zeta_1-\zeta_4|^{\half-\alpha}\right) - \left(
|\zeta_2+\zeta_4|^{\half-\alpha} - |\zeta_4|^{\half-\alpha} \right) =O(\zeta_1 \zeta_2\ \cdot\
|\zeta_4|^{-3/2-\alpha}).\EEQ

Integrating the square of the renormalized amplitude yields (going {\em down} the scales {\em above}
$\zeta_1$) \BEA &&                \zeta_1^2  \left(\int_{|\zeta_1|}^{\infty} \zeta_2^2
 d\zeta_2  \left(\int_{|\zeta_2|}^{\infty}  |\zeta_4|^{-4-4\alpha} d\zeta_4
   \left(\int_{|\zeta_4|}^{\infty}  |\zeta_3|^{-2-4\alpha} d\zeta_3\right)\right)\right)
   \nonumber\\
&&\lesssim \zeta_1^2 \int_{|\zeta_1|}^{\infty} \zeta_2^2
 d\zeta_2  \int_{|\zeta_2|}^{\infty}  |\zeta_4|^{-5-8\alpha} d\zeta_4 \nonumber\\
&&\lesssim   \zeta_1^2 \int_{|\zeta_1|}^{\infty} |\zeta_2|^{-2-8\alpha}  d\zeta_2\nonumber\\ &&=
O(|\zeta_1|^{1-8\alpha}).\EEA Note that the exponents are sufficiently {\em negative} so that these
ultra-violet integrals converge.

 The computation of the integrals yields the same bound as
\BEA &&  M^{j(\zeta_1)\omega^*(G)} \sum_{j(\zeta_2)>j(\zeta_1)} M^{(j(\zeta_2)-j(\zeta_1))
\omega^*(G^{j(\zeta_2)})} \ \cdot \nonumber\\ && \qquad \cdot\ \sum_{j(\zeta_4)>j(\zeta_2)}
M^{(j(\zeta_4)-j(\zeta_2))\omega^*(G^{j(\zeta_4)})}
 \sum_{j(\zeta_3)>j(\zeta_4)}
M^{(j(\zeta_3)-j(\zeta_4))\omega^*(G^{j(\zeta_3)})} , \nonumber\\ \EEA see eq. (\ref{eq:4.5}), since
$\omega^*(G^{j(\zeta_3)})=(1-4\alpha)-2=-1-4\alpha$, $\omega^*(G^{j(\zeta_4)})= -4-8\alpha$ (due to
the fact that the subdiagram with lines $(\xi_4,\zeta_4,\xi_1)$ is renormalized twice),
$\omega^*(G^{j(\zeta_2)})=(1-8\alpha)-2=-1-8\alpha$ and $\omega^*(G)=\omega(G)=1-8\alpha$.

\item Consider now the second Gallavotti-Nicol\`o tree, see Fig. \ref{Fig-renbis4}. One may choose
    as integration variables $L(G)\setminus L'(G)=\{\xi_1,\zeta_2,\zeta_3\}$, so that
    $\zeta_4=\xi_4=\zeta_1-\zeta_2-\xi_1$, $\xi_2=\zeta_2-\zeta_3$, $\xi_3=\zeta_3$. Hence \BEQ
    A(G)=\int d\xi_1 d\zeta_2 d\zeta_3  \left(  |\zeta_1-\zeta_2-\xi_1|^{-\half-\alpha} \ \cdot\
    |\zeta_3|^{-\half-\alpha} |\zeta_2-\zeta_3|^{\half-\alpha} \ \cdot\ |\xi_1|^{\half-\alpha}\
    \cdot\ \zeta_2^{-1} \right)^2.\EEQ

The subdiagram with lines $(\zeta_1,\zeta_4,\xi_4)$ has one external $\phi$-leg, $\xi_1$, hence
needs not be renormalized. On the other hand, the subdiagrams with lines
$(\zeta_1,\zeta_4,\xi_4,\xi_1)$ and $(\xi_2,\zeta_3,\xi_3)$  must be
 renormalized by subtracting their values at $\zeta_2=0$. Hence
  $|\zeta_1-\zeta_2-\xi_1|^{\half-\alpha}$ is replaced with
  $|\zeta_1-\zeta_2-\xi_1|^{-\half-\alpha}-|\zeta_1-\xi_1|^{-\half
-\alpha}=O(\zeta_2\ \cdot |\zeta_1|^{-\frac{3}{2}-\alpha})$, and $|\zeta_2-\zeta_3|^{\half-\alpha}$
by $|\zeta_2-\zeta_3|^{\half-\alpha}-|\zeta_3|^{\half-\alpha}=O(\zeta_2\ \cdot\
|\zeta_3|^{-\half-\alpha}).$

Integrating the square of the renormalized amplitude yields (going {\em up} the scales {\em below}
$\zeta_1$) \BEQ                |\zeta_1|^{-3-2\alpha}  \left(\int_0^{\zeta_1} |\zeta_3|^{-2-4\alpha}
d\zeta_3  \left(\int^{\zeta_3}_0 |\xi_1|^{1-2\alpha} d\xi_1  \left(\int_0^{|\xi_1|} \zeta_2^{2}
d\zeta_2\right)\right)\right)=O(|\zeta_1|^{1-8\alpha}).\EEQ Note that the exponents are sufficiently
{\em positive} so that these infra-red integrals converge.

In order to make the connection with eq. (\ref{eq:4.8}), we replace $(|\zeta_1|^{-3/2-\alpha}
\zeta_2.|\xi_1|^{\half-\alpha})^2=(|\zeta_1|^{-1/2-\alpha}|\xi_1|^{1/2-\alpha}\ \cdot\
\frac{\zeta_2}{\zeta_1})^2$ with $(|\zeta_1|^{-1/2-\alpha} |\xi_1|^{1/2-\alpha}\ \cdot\
\frac{\zeta_2}{\xi_1})^2=|\zeta_1|^{-1-2\alpha}|\xi_1|^{-1-2\alpha} \zeta_2^2.$ The reduced spring
factor $\frac{\zeta_2}{\xi_1}$ takes into account the difference between the {\em minimum scale} of
the diagram with lines $(\zeta_1,\zeta_4,\xi_4,\xi_1)$ and its external leg $\zeta_2$, corresponding
to the lifetime of this diagram; it is the factor which is counted in the multi-scale estimates. The
actual spring factor $\frac{\zeta_2}{\zeta_1}$, which is better, is due to the difference of scales
between the scale where the vertex connecting $\zeta_1,\zeta_4,\xi_1$ and $\zeta_2$ appears and the
scale of the external leg $\zeta_2$. With this slight modification, one gets

\BEQ |\zeta_1|^{-1-2\alpha}  \int_0^{\zeta_1} |\zeta_3|^{-2-4\alpha} d\zeta_3 \int^{\zeta_3}_0
|\xi_1|^{-1-2\alpha} d\xi_1  \int_0^{|\xi_1|} \zeta_2^{2} d\zeta_2.\EEQ

This is equivalent to the  bound given in eq. (\ref{eq:4.8}), \BEA &&
M^{j(\zeta_1)\omega^*(G^{j(\zeta_4)})}  \sum_{j(\zeta_3)<j(\zeta_1)} M^{j(\zeta_3)
 (\omega^*(G^{j(\zeta_3)})-\omega^*(G^{j(\zeta_1)}))} \ \cdot \nonumber\\
&& \quad \sum_{j(\xi_1)<j(\zeta_3)} M^{j(\xi_1) (\omega^*(G^{j(\xi_1)})-\omega^*(G^{j(\zeta_3)})) }
 \sum_{j(\zeta_2)<j(\xi_1)} M^{j(\zeta_2)
(\omega^*(G)-\omega^*(G^{j(\xi_1)})) } , \nonumber\\  \EEA
 since $\omega^*(G^{j(\zeta_4)})=\omega^*(G^{j(\zeta_1)})=-1-2\alpha$, $\omega^*(G^{j(\zeta_3)})=
-6\alpha-2$,  $\omega^*(G^{j(\xi_1)})=-8\alpha-2$ and $\omega^*(G)=\omega(G)=1-8\alpha$.

\end{enumerate}


\section{Proof of H\"older regularity for renormalized skeleton integrals}


We want to prove that, for any indices $(\ell(1),\ldots,\ell(n))$ and $n\le \lfloor 1/\alpha\rfloor$,
\BEQ \Var J^{ts}_B(\ell(1),\ldots, \ell(n))\lesssim |t-s|^{2n\alpha},\EEQ
where $J^{ts}_B$ is defined in Proposition  \ref{lem:barchi-char} and Definition
\ref{def:3.2}.

Consider some multiple contraction $(i_1 i_2),\ldots,(i_{2p-1} i_{2p})$ -- assuming that
$\ell(i_1)=\ell(i_2)$, \ldots, $\ell(i_{2p-1})=\ell(i_{2p})$ --  and the associated contracted integral
$J^{ts}_B(\T;(i_1 i_2), \ldots,(i_{2p-1} i_{2p}))$, where $\T=(\ell(1) \ldots \ell(n))$. By arguments
which may be found in \cite{Unt-fBm}, \S 4.1 (see eq. (4.4) in particular), denoting by $: \ :$ the Wick product of Gaussian variables,
 \BEQ \Var : J^{ts}_B(\T;(i_1 i_2),\ldots,(i_{2p-1}
i_{2p})): \  \le  n!\  \cdot\  \Var J^{ts}_B(\T';(i_1 i_2),\ldots, (i_{2p-1} i_{2p})),\EEQ where $\T'$ has
decorations $(\ell'(1)\ldots \ell'(n))$ such that $\ell'(i)\not=\ell'(j)$ if $i\not=j$ {\em except} if
$\{i,j\}= \{i_{2m-1},i_{2m}\}$ is a pair contraction, as in Lemma \ref{lem:int-Feynman}. Hence, by
Wick's lemma, it suffices to prove that $\Var J^{ts}_B(\T';(i_1 i_2),\ldots, (i_{2p-1} i_{2p}))\lesssim
|t-s|^{2n\alpha}.$

By eq. (\ref{eq:1.35}) and (\ref{eq:1.37}), $J^{ts}_B(\T'; (i_1 i_2),\ldots,(i_{2p-1} i_{2p}))$ is a sum
of terms of the form \BEQ \int d\mathbf{\xi} \left[\prod_{q=1}^p {\cal R}
\Phi^{ts}\right]((\T_q),\mathbf{\xi};(\mathbf{v}_q),(\T'_{q,j});(i_1 i_2), \ldots, (i_{2p-1} i_{2p})),
\label{eq:5:*} \EEQ
 with (following the notations of Lemma \ref{lem:int-Feynman})
\BEQ \left[\prod_{q=1}^p {\cal R}\Phi^{ts}\right](\ .\ )=\left[ \del {\cal R} \overline{\SkI}^{ts} \
\cdot\ {\cal R} \overline{\SkI}^s \right]_{\hat{\nu}(\mathbf{\xi})} \left( \prod_{q=1}^p
\Roo_{\mathbf{v}_q}\T_q,\prod_{q=1}^p \prod_j \T'_{q,j}; (i_1 i_2),\ldots, (i_{2p-1} i_{2p}) \right).\EEQ
The contractions induce links between some of the trees $(\T_q)_q$, $(\T'_{q,j})_{q,j}$. The resulting
 connected components
may be represented by Feynman graphs of two types: (i) {\em ``rooted'' Feynman diagrams}
$G_1,\ldots,G_I$
 containing
some (possibly many) root part $\Roo_{\mathbf{v}_q}\T_q$; (ii) {\em ``unrooted Feynman diagrams}
$G'_1,\ldots,G'_{I'}$ containing only leaf parts of type $\T'_{q,j}$. It turns out that the unconvenient
vertex-decorating characteristic function ${\bf 1}_{|\xi_1|\le \ldots\le |\xi_n|}$ may be replaced with
the following much simpler characteristic function $f$. Let $G_1$ be the rooted diagram
 containing $\xi_1$. For every unrooted diagram $G'_i$ of type (ii), choose some $\xi$-leg
$\xi'_i$ belonging to $G'_i$ and let $f_i:={\bf 1}_{j(\xi'_i)\ge j(\xi_1)}$. Then set
$f=f^{j(\xi_1)}:=\prod_{i=1}^{I'} f_i(\mathbf{\xi})$. The integral $\int d\mathbf{\xi} \ {\bf 1}_{|\xi_1|\le
\ldots\le |\xi_n|} (\ .\ )$
 in (\ref{eq:1.35}) is now
 replaced by a simple
sum $\sum_{j=-\infty}^{+\infty} (\ .\ )$, with $j=j(\xi_1)$, and $\hat{\nu}(\mathbf{\xi})$ by a measure
depending only on the scale $j(\xi_1)$, \BEA  \nu^j &:=&\sum_{j'_i\ge j,i=1,\ldots,I'} \int d\xi_1\ldots
d\xi_n {\bf 1}_{j(\xi_1)=j} \left[ \prod_{i=1}^{I'} {\bf 1}_{j(\xi'_i)=j'_i} \right] \otimes_{k=1}^n
{\cal F}(\Gamma'(\ell\circ\sigma(k)))(\xi_k) \nonumber\\ &&= \int_{M^j\le |\xi_1|\le M^{j+1}} d\xi_1
\int d\xi_2\ldots d\xi_n \left[\prod_{i=1}^{I'} {\bf 1}_{|\xi'_i|\ge M^j} \right] \otimes_{k=1}^n {\cal
F}(\Gamma'(\ell\circ\sigma(k)))(\xi_k). \nonumber\\ \EEA Since $f(\mathbf{\xi})\le {\bf 1}_{|\xi_1|\ge
\ldots\ge |\xi_n|}$, the associated renormalized quantity \BEQ \sum_{j=-\infty}^{+\infty}
\left[\prod_{q=1}^p {\cal R}\Phi^{ts} \right]((\T_q),j;(\mathbf{v}_q),
 (\T'_{q,j}); (i_1 i_2),\ldots,(i_{2p-1} i_{2p})) \EEQ
has a larger variance than the original one, eq. (\ref{eq:5:*}),  contributing to $J^{ts}_B$.

\bigskip

The purpose of this section is to prove the estimates
\BEQ \Var\left( \sum_{j=-\infty}^{+\infty}
\left[\prod_{q=1}^p {\cal R}\Phi^{ts} \right]((\T_q),j;(\mathbf{v}_q), (\T'_{q,j}); (i_1
i_2),\ldots,(i_{2p-1} i_{2p})) \right) \lesssim |t-s|^{2n\alpha}, \label{eq:1791}
\EEQ
 from which Theorem \ref{th:0.1}
follows. They are a simple consequence of Theorem \ref{th:4.1} and of the following two lemmas.

\begin{Lemma} [bound for bilateral ``rooted'' diagrams]  \label{lem:5.1} (see Lemma \ref{lem:int-Feynman} (2) for
notations)

Let $q\ge 1$ and $q'\ge 0$, $n:=q+q'$, and $n'\ge 0$ such that $n+n'< \lfloor 1/\alpha \rfloor$. Rename
$(\zeta_{r_1},\ldots,\zeta_{r_q},\zeta_{r'_1},\ldots,\zeta_{r'_{q'}})$,
 resp. $(\bar{\zeta}_{r_1},\ldots,\bar{\zeta}_{r_q},\bar{\zeta}_{r'_1},\ldots,\bar{\zeta}_{r'_{q'}})$,
 as $\zeta_1,\ldots,\zeta_n$, resp. $\tilde{\zeta}_1,\ldots,\tilde{\zeta}_n$, so that
 $|\zeta_1|<\ldots<|\zeta_n|$ and
$|\tilde{\zeta}_1|<\ldots<|\tilde{\zeta}_n|$, and let $\zeta_{ext}:=\sum_{m=1}^n \zeta_m$,
$\tilde{\zeta}_{ext}:=\sum_{m=1}^n \tilde{\zeta}_m$.

 Let \BEA && I(j_{ref}):=\int_{|\xi_{ext}|\le
M^{j_{ref}}} d\xi_{ext} \int d\mathbf{\zeta}_{ext} \int d\tilde{\mathbf{\zeta}}_{ext}
 \del(\zeta_{ext}=\xi_{ext})\del(\tilde{\zeta}_{ext}=-\xi_{ext}) \nonumber\\
 && \qquad  M^{(1-2(n+n')\alpha)j_{ref}} \left(
 \frac{\min(|\zeta_1|,M^{j_{ref}})}{\max(|\zeta_n|,M^{j_{ref}})}
 \right)^{\alpha^-/2}  \left(
 \frac{\min(|\tilde{\zeta}_1|,M^{j_{ref}})}{\max(|\tilde{\zeta}_n|,M^{j_{ref}})}
 \right)^{\alpha^-/2}
\nonumber\\ && \left( \prod_{m=1}^q \left| \frac{e^{\II t\zeta_{r_m}}-e^{\II
s\zeta_{r_m}}}{\zeta_{r_m}}\right| \prod_{m'=1}^{q'} \left| \frac{1}{\zeta_{r'_{m'}}} \right| \right)
\left( \prod_{m=1}^q \left| \frac{e^{\II t\bar{\zeta}_{r_m}}-e^{\II
s\bar{\zeta}_{r_m}}}{\bar{\zeta}_{r_m}}\right| \prod_{m'=1}^{q'} \left| \frac{1}{\bar{\zeta}_{r'_{m'}}}
\right| \right). \nonumber\\ \EEA

Then \BEQ \sum_{j_{ref}=-\infty}^{+\infty} I(j_{ref})\lesssim |t-s|^{2(n+n')\alpha}.   \EEQ

\end{Lemma}

{\bf Proof.}

\begin{itemize}

\item[(i)] $(M^{j_{ref}}<\frac{1}{|t-s|})$

 Integrate first over the variables
{\em larger} than $\frac{1}{|t-s|}$  -- which defines the {\em ultra-violet range}
in this situation --, say $|\zeta_n|>\ldots>|\zeta_{k+1}|$ and
$|\tilde{\zeta}_n|>\ldots>|\tilde{\zeta}_{\tilde{k}+1}|$. Let for instance $|{\zeta}_n|>|\tilde{\zeta}_n|$.
The integral $\int_{|\zeta_n|>|\zeta_{n-1}|} \frac{d\zeta_n}{|\zeta_n|} {\bf 1}_{|\zeta_{ext}|\le
M^{j_{ref}}}$ is bounded up to a constant by $\frac{M^{j_{ref}}}{|\zeta_{n-1}|}$, and (due to momentum conservation)
the integral over $\tilde{\zeta}_n$ is not performed. Then

\BEQ \int_{|\tilde{\zeta}_{\tilde{k}+1}|>\frac{1}{|t-s|}}
\frac{d\tilde{\zeta}_{\tilde{k}+1}}{|\tilde{\zeta}_{\tilde{k}+1}|} \ldots
\int_{|\tilde{\zeta}_{n-1}|>|\tilde{\zeta}_{n-2}|}
 \frac{d\tilde{\zeta}_{n-1}}{|\tilde{\zeta}_{n-1}|} \
 \cdot\ |\tilde{\zeta}_{n-1}|^{-1} =O(|t-s|) \label{eq:5.22}\EEQ

and similarly for the untilded integrals, with an extra $M^{j_{ref}}$ factor.

Integrating in the {\em infra-red range}, namely, over the variables smaller than $\frac{1}{|t-s|}$ (if any) yields then, using the $\alpha^-$-spring factors,
\BEQ
|t-s|^{\alpha^-/2} \int_{|\tilde{\zeta}_{\tilde{k}}|<\frac{1}{|t-s|}}
\frac{d\tilde{\zeta}_{\tilde{k}}}{|\zeta_{\tilde{k}}|} \ldots
\int_{|\tilde{\zeta}_2|<|\tilde{\zeta}_3|} \frac{d\tilde{\zeta}_2}{|\tilde{\zeta}_2|}
\int_{|\tilde{\zeta}_1|<|\tilde{\zeta}_2|} \frac{d\tilde{\zeta}_1}{|\tilde{\zeta}_1|}\ \cdot\
|\tilde{\zeta}_1|^{\alpha^-/2}=O(1)  \EEQ
and similarly for the untilded integrals.

\bigskip

The above arguments do not hold if {\em all} variables are
smaller than $\frac{1}{|t-s|}$.
 Then  one must use the hypothesis
that at least one of the $\zeta$-variables, say, $\zeta_k$, is accompanied by the factor
$|\frac{e^{\II t\zeta_k}-e^{\II s\zeta_k}}{\zeta_k}|=O(|t-s|)$ instead of $O(\frac{1}{|\zeta_k|})$,
and similarly for some $\tilde{\zeta}$-variable, say, $\tilde{\zeta}_{\tilde{k}}$. One computes
 \BEQ
\int_{|\tilde{\zeta}_{\tilde{k}-1}|<|\tilde{\zeta}_{\tilde{k}}|}
\frac{d\tilde{\zeta}_{\tilde{k}-1}}{|\zeta_{\tilde{k}-1}|} \ldots
\int_{|\tilde{\zeta}_2|<|\tilde{\zeta}_3|} \frac{d\tilde{\zeta}_2}{|\tilde{\zeta}_2|}
\int_{|\tilde{\zeta}_1|<|\tilde{\zeta}_2|} \frac{d\tilde{\zeta}_1}{|\tilde{\zeta}_1|}\ \cdot\
|\tilde{\zeta}_1|^{\alpha^-/2}=O(|\tilde{\zeta}_{\tilde{k}}|^{\alpha^-/2}) \EEQ and
similarly for the untilded integrals, and
 \BEA &&
 |\tilde{\zeta}_n|^{-1-\alpha^-/2} \int_{|\tilde{\zeta}_{n-1}|<|\tilde{\zeta}_n|} \frac{d\tilde{\zeta}_{n-1}}{|\tilde{\zeta}_{n-1}|}
\int_{|\tilde{\zeta}_{n-2}|<|\tilde{\zeta}_{n-1}|}
 \frac{d\tilde{\zeta}_{n-2}}{|\tilde{\zeta}_{n-2}|} \ldots \nonumber\\
 && \qquad
 \int_{|\tilde{\zeta}_{\tilde{k}}|<|\tilde{\zeta}_{\tilde{k}+1}|}  d\tilde{\zeta}_{\tilde{k}}
 |t-s|
\ \cdot\  |\tilde{\zeta}_{\tilde{k}}|^{\alpha^-/2}=O(|t-s|),  \nonumber\\ \EEA

 \BEA &&
  \int_{|{\zeta}_{n}|<\frac{1}{|t-s|}} \frac{d{\zeta}_{n}}{|\zeta_{n}|^{1+\alpha^-/2}}
  {\bf 1}_{|\zeta_{ext}|\le M^{j_{ref}}}
\int_{|{\zeta}_{n-1}|<|{\zeta}_{n}|}
 \frac{d{\zeta}_{n-1}}{|{\zeta}_{n-1}|} \ldots \nonumber\\ &&\qquad
 \int_{|{\zeta}_{k}|<|{\zeta}_{k+1}|}  d\zeta_{k}
 |t-s|
\ \cdot\  |{\zeta}_{k}|^{\alpha^-/2}=O( M^{j_{ref}} |t-s|).  \nonumber\\ \EEA

All together (in both cases) : $I(j_{ref})\lesssim M^{(1-2(n+n')\alpha)j_{ref}} \ \cdot\ |t-s|\ \cdot\ (M^{j_{ref}}|t-s|)=
M^{(2-2(n+n')\alpha)j_{ref}} |t-s|^2.$  Since by assumption $M^{j_{ref}}<\frac{1}{|t-s|}$, this sums up
to $\sum_{j_{ref}<\log_M \frac{1}{|t-s|}} I(j_{ref})\lesssim |t-s|^{2(n+n')\alpha}.$

\item[(ii)]  $(M^{j_{ref}}>\frac{1}{|t-s|})$

The arguments of (i) may be repeated word for word, except that the ultra-violet range is now defined by $|\zeta|,|\tilde{\zeta}|>M^{j_{ref}}$. Then $I(j_{ref})\lesssim M^{-2(n+n')\alpha j_{ref}}$ and $\sum_{j_{ref}>\log_M \frac{1}{|t-s|}} M^{-2(n+n')\alpha j_{ref}}=O(|t-s|^{2(n+n')\alpha}).$

\end{itemize}

The mixed cases, when e.g $\frac{1}{|t-s|}$ is large with respect to the $\zeta$-variables but small
with respect to the $\tilde{\zeta}$-variables, are treated in the same way and left to the reader.

 \hfill \eop

\begin{Lemma}[bound for bilateral ``unrooted'' diagrams]  \label{lem:5.2}
(same notations as in Lemma
\ref{lem:5.1}).
 Assume $q=0$, so that $n':=q'<\lfloor 1/\alpha \rfloor$. Then
\BEQ I(j_{ref})
\lesssim M^{-2n'\alpha j_{ref}}.  \EEQ

\end{Lemma}

Lemma \ref{lem:5.2} has already been proved, as part of Lemma \ref{lem:5.1} (ii).

These two lemmas extend with very minor changes to unilateral diagrams.

\vskip 2cm

We may now easily finish the proof of the estimates eq. (\ref{eq:1791}). Lemma \ref{lem:5.2} yields an estimate for renormalized skeleton
integrals associated to ``unrooted'' diagrams $G'_i$, $i=1,\ldots,I'$, where some reference scale $j'_i=j(\xi_{ref})$
has been chosen according to the rules of Theorem \ref{th:4.1}. Summing over all scales $j'_i\ge j(\xi_1)$ -- where $|\xi_1|$ is the smallest $\xi$-variable, as in the introduction to the present section -- yields $O(M^{-2n'\alpha j(\xi_1)})$. Then
the product of  factors $M^{-2n'\alpha j(\xi_1)}=\prod_{i=1}^{I'} M^{-2n'_i \alpha j(\xi_1)}$ associated to all unrooted
diagrams $(G'_i)_{i=1,\ldots,I'}$ is "grafted" into the rooted diagram $G_1$ containing $\xi_1$.

\medskip

Turn now to the rooted diagrams $G_1,\ldots,G_I$. Choose $j(\xi_1)$ as reference scale for $G_1$ if $G_1$ is unilateral; choose some reference scale according to the rules of Theorem \ref{th:4.1} for $G_2,\ldots,G_I$, and for $G_1$ if $G_1$ is bilateral. Then apply Lemma \ref{lem:5.2}. \hfill\eop

\vskip 5 cm

Let us add two comments to finish with.

\begin{enumerate}

\item The most "tricky" part in the story is obviously the infra-red behaviour of Feynman diagrams, particularly when $n$ is {\em large}, i.e. $n>\frac{1}{2\alpha}$. The infra-red convergence of these "large" diagrams is ensured by the somewhat complicated interplay between half-diagrams and full diagrams, the key point being the existence of small enough spring factors. In a previous attempt, we tried to use the BPHZ renormalization scheme associated to the Connes-Kreimer algebra $\bf H$, instead of considering the associated Feynman half-diagrams. The coproduct of $\bf H$ is much simpler than that of Feynman diagrams. Unfortunately, some "large" diagrams are infra-red divergent.

\item The results of this article may probably be extended to an arbitrary $\alpha$-H\"older path $\Gamma$, by rewriting $\Gamma$ as $I_{\alpha^-}(D_{\alpha^-}(\Gamma))$, where $I_{\alpha^-}$, resp. $D_{\alpha^-}$ are fractional
    integration, resp. derivation operators, and $\alpha^-<\alpha$. Then what one should really do is renormalize iterated fractional integration operators, while $\Gamma$ would only play a "decorative" r\^ole; see Remark 1. after Definition \ref{def:3.2}. The construction would  make use of Besov norms as in \cite{Unt-Holder}.

\end{enumerate}

\vskip 2cm

{\bf Acknowledgements.} We wish to thank Lo\"\i c Foissy, Kurusch Ebrahimi-Fard and Dominique
 Manchon for useful discussions on the subject and for enlightenments on the algebraic part of the construction.



\newpage


\begin{thebibliography}{99}


\bibitem{BHL}  R. F. Bass, B. M. Hambly, T. J. Lyons. {\it Extending the Wong-Zakai theorem to
    reversible Markov processes}, J. Eur. Math. Soc. {\bf 4}, 237--269 (2002).




\bibitem{BroFra03} C. Brouder, A. Frabetti. {\it QED Hopf algebras on planar binary trees}, Journal of
    Algebra {\bf 267}, 298--322 (2003).

\bibitem{BroFraKra06} C. Brouder, A. Frabetti, C. Krattenthaler, {\it Non-commutative Hopf algebra of
    formal diffeomorphisms}, Advances in Math. {\bf 200}, 479--524 (2006).

\bibitem{But72} J. C. Butcher. {\it An algebraic theory of integration methods}, Math. Comp. {\bf 26},
    79--106 (1972).

\bibitem{CEM08} D. Calaque, K. Ebrahimi-Fard, D. Manchon. {\it Two Hopf algebras of trees interacting}.
    Preprint arXiv:0806.2238.

\bibitem{ChaLiv07} F. Chapoton, M. Livernet. {\em Relating two Hopf algebras built from an operad},
    International Mathematics Research Notices, Vol. 2007, Article ID rnm131.



\bibitem{ConKre98} A. Connes, D. Kreimer. {\it Hopf algebras, renormalization and non-commutative
    geometry}, Comm. Math. Phys. {\bf 199} (1), 203--242 (1998).

\bibitem{ConKre00} A. Connes, D. Kreimer. {\it Renormalization in quantum field theory and the
    Riemann-Hilbert problem (I)}, Comm. Math. Phys. {\bf 210} (1), 249--273 (2000).

\bibitem{ConKre01}  A. Connes, D. Kreimer. {\it Renormalization in quantum field theory and the
    Riemann-Hilbert problem (II)}, Comm. Math. Phys. {\bf 216} (1), 215--241 (2001).

\bibitem{CQ02} L. Coutin, Z. Qian. {\it Stochastic analysis, rough path analysis and fractional Brownian
    motions}, Probab. Theory Related Fields  {\bf 122}  (2002),  no. 1, 108--140.


\bibitem{Foi02} L. Foissy. {\it Les alg\`ebres de Hopf des arbres enracin\'es d\'ecor\'es (I)}, Bull.
    Sci. Math., {\bf 126} (3), 193--239, {\it and (II)}, Bull. Sci. Math., {\bf 126} (4), 249--288
    (2002).

\bibitem{FoiUnt} L. Foissy, J. Unterberger. {\it Ordered forests, permutations and iterated integrals}.
    Preprint arXiv:1004.5208 (2010).

\bibitem{FV} P. Friz, N. Victoir: {\it Multidimensional dimensional processes seen as rough paths}.
    Cambridge University Press (2010).



\bibitem{Gar} A. Garsia. {\it Continuity properties of Gaussian processes with multidimensional time
parameter}, Proceedings of the Sixth Berkeley Symposium on Mathematical Statistics and Probability
Vol. II: Probability theory, 369-374. Univ. California Press (1972).

\bibitem{GNRV} M. Gradinaru, I. Nourdin, F. Russo, P. Vallois. {\em $m$-order integrals and generalized
    It\^o's formula: The case of a fractional Brownian motion with any Hurst index}, Ann. Inst. H.
    Poincar\'e Probab. Statist. {\bf 41}, 781--806 (2005).


\bibitem{Gu} M. Gubinelli. {\it Controlling rough paths}, J. Funct. Anal. {\bf 216}, 86--140 (2004).

\bibitem{Gu2} M. Gubinelli. {\it Ramification of rough paths}. Preprint available on {\it Arxiv}
    (2006).



\bibitem{Hepp} K. Hepp. {\it Proof of the Bogoliubov-Parasiuk theorem on renormalization}, Commun. Math.
    Phys. {\bf 2} (4), 301--326 (1966).



\bibitem{HL} B. Hambly, T. J. Lyons. {\it Stochastic area for Brownian motion on the Sierpinski basket},
    Ann. Prob. {\bf 26} (1), 132--148 (1998).


\bibitem{Kre99}  D. Kreimer. {\it Chen's iterated integral represents the operator product expansion},
    Adv. Theor. Math. Phys. {\bf 3} (3), 627--670 (1999).

\bibitem{Lej03} A. Lejay. {\it An introduction to rough paths}, S\'eminaire de probabilit\'es XXXVII,
    Lecture Notes in Mathematics, Springer (2003).



\bibitem{LQ-bk} T. Lyons, Z. Qian (2002): {\it System control and rough paths}. Oxford University Press
    (2002).

\bibitem{LyoVic07} T. Lyons, N. Victoir. {\it An extension theorem to rough paths}, Ann. Inst. H.
    Poincar\'e Anal. Non Lin\'eaire {\bf 24} (5), 835--847 (2007).


\bibitem{MagUnt} J. Magnen, J. Unterberger. {\em From constructive field theory to fractional stochastic
    calculus. (I) The L\'evy area of fractional Brownian motion with Hurst index
    $\alpha\in(\frac{1}{8},\frac{1}{4})$}. Preprint arXiv:1004.5208.

\bibitem{Mur} A. Murua. {\it The shuffle Hopf algebra and the commutative Hopf algebra of labelled
    rooted trees.} Available on www.ehu.es/ccwmuura/research/shart1bb.pdf.

\bibitem{Mur2} A. Murua. {\it The Hopf algebra of rooted trees, free Lie algebras, and Lie series},
    Found. Comput. Math. {\bf 6} (4), 387--426 (2006).


\bibitem{Nu-cours} D. Nualart. {\it Stochastic calculus with respect to the fractional Brownian motion
    and applications}, Contemporary Mathematics {\bf 336}, 3-39 (2003).




\bibitem{Riv} V. Rivasseau. {\it From Perturbative to Constructive Renormalization}, Princeton Series in
    Physics (1991).




\bibitem{TinUnt08} S. Tindel, J. Unterberger. {\it The rough path associated to the multidimensional
    analytic fBm with any Hurst parameter}. Preprint available on {\it Arxiv} (2008). A
    para\^\i tre dans: Collectanea Mathematica.









\bibitem{Unt-1} J. Unterberger. {\it Stochastic calculus for fractional Brownian motion with Hurst
    parameter $H>1/4$: a rough path method by analytic extension}, Ann. Prob. {\bf 37} (2), 565--614
    (2009).



\bibitem{Unt-fBm} J. Unterberger. {\it A rough path over  multi-dimensional fractional Brownian motion
    with arbitrary Hurst index by Fourier normal ordering}, Stoch.  Proc. Appl. {\bf 120} (8), 1444-1472 (2010).  

\bibitem{Unt09ter} J. Unterberger. {\it A L\'evy area by Fourier normal ordering for multidimensional
    fractional Brownian motion with small Hurst index}.  Preprint arXiv:0906.1416.


\bibitem{Unt-Holder} J. Unterberger. {\it H\"older-continuous rough paths by Fourier normal ordering}, Comm.
Math. Phys. {\bf 298} (1), 1--36 (2010).

\bibitem{Unt-lin} J. Unterberger. {\it Moment estimates for solutions of linear stochastic differential
equations driven by analytic fractional Brownian motion}. Preprint arXiv:0905.0782.
   

\bibitem{FVT} F. Vignes-Tourneret. {\em Renormalisation des th\'eories de
champs non commutatives},Th\`ese de doctorat de l'Universit\'e Paris 11,
arXiv:math-ph/0612014.

\bibitem{Wal00} M. Waldschmidt. {\it Valeurs z\^eta multiples. Une introduction}, Journal de Th\'eorie
    des Nombres de Bordeaux {\bf 12} (2), 581--595 (2000).

\end{thebibliography}
\end{document}